\documentclass{article}
\usepackage{amssymb}
\usepackage{amsmath}
\usepackage{amscd}
\usepackage{amsgen}
\usepackage{amsopn}
\usepackage{amsthm}
\usepackage{eucal}
\usepackage{graphicx}
\usepackage{amsfonts,latexsym}

\usepackage[all]{xy}           

\usepackage{xspace}
\usepackage{epsfig}
\usepackage{float}

\usepackage{verbatim}
\usepackage{upref}

\theoremstyle{plain}
\newtheorem{thm}{Theorem}[section]
\newtheorem{lem}[thm]{Lemma}
\newtheorem{prop}[thm]{Proposition}
\newtheorem{cor}[thm]{Corollary}

\theoremstyle{definition}

\newtheorem{defn}[thm]{Definition}
\newtheorem{rem}[thm]{Remark}

\newtheorem{eg}[thm]{Example}

\DeclareMathOperator{\Hom}{{Hom}}

\font\cyr=wncyr10

\newcommand{\nc}{\newcommand}

 \def\wvec#1#2{\Big[ {\displaystyle #1 \atop \displaystyle #2}\Big]}
 \nc{\subv}{{^{\star}}}
 \nc{\msh}{{\ast_q}}
 \nc{\lms}{{\longmapsto}}
 \nc{\lra}{{\longrightarrow}}
 \nc{\hcirc}{\hat{\circ}}
 \nc{\id}{{\rm{id}}}
 \nc{\im}{\rm{im}}
 \nc{\incl}{\rm{incl}}
 \nc{\sha}{{\mbox{\cyr x}}}
 \nc{\us}{\underset}
 \nc{\tn}{{\tilde{n}}}
 \nc{\ra}{\rightarrow}
 \nc{\vs}{{\vec{s}}}
 \nc{\va}{{\vec{a}}}
 \nc{\vb}{{\vec{b}}}
 \nc{\vc}{{\vec{c}}}
 \nc{\vf}{{\vec{f}}}
 \nc{\vd}{{\vec{\Delta}}}
 \nc{\vx}{{\vec{x}}}
 \nc{\vy}{{\vec{y}}}
\nc{\vi}{{\vec{\imath}}}
 \nc{\vone}{{\vec{1}}}
 \nc{\vo}{{\vec{0}}}
 \nc{\vr}{{\vec{r}}}
 \nc{\os}{{\overset}}
 \nc{\Z}{{\mathbb Z}}
 \nc{\R}{{\mathbb R}}
 \nc{\N}{{\mathbb N}}
 \nc{\ZN}{{\mathbb Z_{\ge 0}}}
 \nc{\Q}{{\mathbb Q}}
 \nc{\C}{{\mathbb C}}
 \nc{\D}{{\mathcal D}}
 \nc{\Om}{{\Omega}}
 \nc{\caT}{{\mathcal T}}
 \nc{\tB}{{\tilde B}}

 \nc{\Li}{{\rm Li}}
 \nc{\wh}{\widehat}
 \nc{\shaq}{{\sha_q\,}}
 \nc{\gam}{{\gamma}}
 \nc{\ga}{{\alpha}}
 \nc{\etaq}{{\eta_q}}
 \nc{\vep}{{\varepsilon}}
 \nc{\gl}{{\lambda}}
 \nc{\gb}{{\beta}}
 \nc{\gd}{{\delta}}
 \nc{\gs}{{\sigma}}
 \nc{\gS}{{\Sigma}}
 \nc{\gz}{{\zeta}}
 \nc{\Zq}{{Z_q}}
 \nc{\uZq}{{{\tilde{Z}_q}}}
 \nc{\gzq}{{\bar \zeta_q}}
 \nc{\bzq}{{{\bar \zeta}_q}}
 \nc{\xiq}{{\xi_q}}
 \nc{\sif}{{\mathcal S}}
 \nc{\gt}{{\tau}}
 \nc{\gk}{{\kappa}}
 \nc{\bw}{\bar w}
 \nc{\Lra}{\Longrightarrow}
 \nc{\fS}{{\mathfrak S}}
 \nc{\DD}{{\mathfrak D}}
 \nc{\Llra}{\Longleftrightarrow}
 \nc{\zq}{{\zeta_q}}
 \nc{\qup}{{q\uparrow 1}}
 \nc{\bfone}{{\bf 1}}
 \nc{\lp}{\Big ( } \nc{\llp}{\Big (} \nc{\Llp}{\left (}
 \nc{\rp}{\Big ) } \nc{\rrp}{\Big )} \nc{\Rrp}{\right )}
 \nc{\lb}{\big < } \nc{\llb}{\!\Big \langle } \nc{\Llb}{\! \left <}
 \nc{\rb}{\big >  } \nc{\rrb}{\Big \rangle \!} \nc{\Rb}{\Big
 \rangle\! } \nc{\length}{{\rm leng}} \nc{\cop}{{\rm cop}}

 \nc{\Alg}{\mathbf{Alg}}
 \nc{\Bax}{\mathbf{Bax}}
 \nc{\bfk}{{\bf k}}
 \nc{\bfe}{{\bf e}}
 \nc{\base}[1]{\bfone^{\otimes ({#1}+1)}} 
 \nc{\Rings}{\mathbf{Rings}}
 \nc{\Sets}{\mathbf{Sets}}

 \nc{\BA}{{\Bbb A}}
 \nc{\EE}{{\Bbb E}}
 \nc{\GG}{{\Bbb G}}
 \nc{\HH}{{\Bbb H}}
 \nc{\LL}{{\Bbb L}}
 \nc{\TT}{{\Bbb T}}
 \nc{\VV}{{\Bbb V}}

 \nc{\cala}{{\mathcal A}}
 \nc{\calc}{{\mathcal C}}
 \nc{\cald}{\mathcal{D}}
 \nc{\cale}{{\mathcal E}}
 \nc{\calf}{{\mathcal F}}
 \nc{\calg}{{\mathcal G}}
 \nc{\calh}{{\mathcal H}}
 \nc{\cali}{{\mathcal I}}
 \nc{\call}{{\mathcal L}}
 \nc{\calm}{{\mathcal M}}
 \nc{\caln}{{\mathcal N}}
 \nc{\calo}{{\mathcal O}}
 \nc{\calp}{{\mathcal P}}
 \nc{\calr}{{\mathcal R}}
 \nc{\calt}{{\mathcal T}}
 \nc{\calw}{{\mathcal W}}
 \nc{\calx}{{\mathcal X}}
 \nc{\CA}{\mathcal{A}}
 \nc{\calz}{{\mathcal Z}}
 \nc{\calL}{{\mathcal L}}
 \nc{\fraka}{{\frak a}}
 \nc{\frakA}{{\frak A}}
 \nc{\frakb}{{\frak b}}
 \nc{\frakB}{{\frak B}}
 \nc{\frakH}{{\frak H}}
 \nc{\frakM}{{\frak M}}
 \nc{\frakS}{{\frak S}}
 \nc{\bfrakM}{\overline{\frakM}}
 \nc{\frakm}{{\frak m}}
 \nc{\frakP}{{\frak P}}
 \nc{\frakN}{{\mathfrak N}}
 \nc{\frakp}{{\frak p}}

\begin{document}

\title{Renormalization of Multiple $q$-Zeta Values}
\author{Jianqiang Zhao}
\date{}
\maketitle
\begin{center}
{\large Department of Mathematics, Eckerd College, St. Petersburg,
FL 33711}
\end{center}

\medskip

\noindent{\bf Abstract.} In this paper we shall define the
renormalization of the multiple $q$-zeta values (M$q$ZV) which are
special values of multiple $q$-zeta functions $\zq(s_1,\dots,s_d)$
when the arguments are all positive integers or all non-positive
integers. This generalizes the work of Guo and Zhang \cite{GuZ} on
the renormalization of Euler-Zagier multiple zeta values. We show
that our renormalization process produces the same values if the
M$q$ZVs are well-defined originally and that these
renormalizations of M$q$ZV satisfy the $q$-stuffle relations if we
use shifted-renormalizations for all divergent
$\zq(s_1,\dots,s_d)$ (i.e., $s_1\le 1$). Moreover, when $\qup$ our
renormalizations agree with those of Guo and Zhang.


\medskip


\section{Introduction}
The Euler-Zagier multiple zeta functions are defined as nested
generalizations of the Riemann zeta function:
\begin{equation}\label{equ:oldzeta}
\zeta(s_1,\dots, s_d):=\sum_{k_1>\dots>k_d>0} k_1^{-s_1}\cdots
k_d^{-s_d}
\end{equation}
for complex variables $s_1,\dots, s_d$ satisfying
$\Re(s_1+\dots+s_j)>j$ for all $j=1,\dots,d$. The special values
of this function at positive integers are called \emph{multiple
zeta values} (MZVs) and were first studied systematically by Euler
\cite{LE1} when $d=2$. Nevertheless, only in the past fifteen
years these values have been found to have significant arithmetic,
algebraic, geometric and physics meanings and have since been
under intensive investigation (see \cite{BK,MG,H1,LM,Zag}).

In another direction in \cite{Zana} we show by using generalized
functions that multiple zeta functions can be analytically
continued to $\C^d$ as a meromorphic function with simple poles.
We will henceforth always refer to this analytic continuation when
we speak of multiple zeta functions in the rest of this paper. The
precise location of the simple poles form the following set (see
\cite{AET}):
\begin{equation}\label{fSd}
\fS_d=\left\{(s_1,\dots,s_d)\in \C^d\left| \aligned
s_1&= 1;, \text{ or } s_1+s_2\in \{1\}\cup 2\Z_{\le 1},\text{ or}\\
s_1&+\dots+s_j\in\Z_{\le j} \text{ for }3\le j\le d
\endaligned
\right.\right\}.
\end{equation}
Hence MZVs at non-positive integers are not always defined.
Recently, Guo and his collaborators (\cite{EG2,GuZ}) have applied
the Rota-Baxter algebra technique to the study of MZVs after
noticing that the stuffle (stuffing+shuffle) relations reflect
exactly the Rota-Baxter property. In \cite{GuZ} the
renormalization is carried out for the MZVs and they show that
when $\zeta(s_1,\dots,s_d)$ is defined then its renormalization
agrees with the value itself, provided that $s_i$'s are all
positive or all non-positive. Moreover, these renormalizations
satisfy the stuffle (or quasi-shuffle) relations. The importance
of this result is related to the conjecture (see \cite{IKZ}) that
to obtain all the relations among MZVs of the same weight it
suffices to use all the double shuffle relations including those
of the renormalization of MZVs at positive integers.

On the other hand, we can define the $q$-analog ($0<q<1$) of
multiple zeta functions as follows (see \cite{Zqmv}). For complex
variables $s_1,\dots, s_d$ satisfying $\Re(s_1+\dots+s_j)>j$ for
all $j=1,\dots,d$, set
\begin{equation}\label{equ:qzeta}
\zq(s_1,\dots,s_d):=\sum_{k_1>\dots>k_d>0}
\frac{q^{k_1(s_1-1)+\cdots+k_d(s_d-1)}}{[k_1]^{s_1}\cdots[k_d]^{s_d}}
\end{equation}
where for any real number $r$ we write $[r]=(1-q^r)/(1-q).$ When
$d=1$ this is the same as the $q$-analog of the Riemann zeta
function defined in \cite{KKW}. By using Euler-Maclaurin
summations we obtained their meromorphic continuations to $\C^d$
with following singularities (which are all simple):
\begin{equation}\label{fSd'}
\fS_d'=\left\{(s_1,\dots,s_d)\in \C^d\left| \aligned &s_1\in
1+\frac{2\pi i}{\ln q}\Z, \text{ or }
s_1\in\Z_{\le 0}+\frac{2\pi i}{\ln q}\Z_{\ne 0}, \\
&\text{or } s_1+\dots+s_j\in \Z_{\le j} +\frac{2\pi i}{\ln q}\Z
\text{ for } j>1
\endaligned \right.\right\}\supset \fS_d.
\end{equation}
Here the last part in $\fS_d'$ is vacuous if $d=1$.  One can see
that these $q$-analogues have much more poles than their ordinary
counterparts. But when $q$ approaches 1 we indeed recover exactly
the poles of the multiple zeta functions. In fact, by \cite[Main
Theorem]{Zqmv} for all $(s_1,\dots,s_d)\in \C^d\setminus \fS_d$
$\lim_\qup \zq(s_1,\dots,s_d)=\zeta(s_1,\dots,s_d)$, which shows
that our $q$-analogue is the correct choice.

The analytical continuation of multiple zeta functions \cite{Zana}
utilizes generalized functions. But Euler-Maclaurin summation can
also be used instead which actually provides the main idea of
special value computations contained in this paper. For future
reference we define the Bernoulli polynomials $B_k(x)$ and its
periodic analogue $\tB_k(x)$ by
\begin{equation}\label{equ:bernoulli}
 \frac{te^{xt}}{e^t-1}=\sum_{n=0}^\infty
 B_n(x)\frac{t^n}{n!},\qquad \tB_k(x)=B_k(\{x\}), x\ge 1,
\end{equation}
where $\{x\}$ is the fractional part of $x$. We can then prove
the analytic continuation of $\zeta(s_1,\dots, s_d)$ using these
functions. See \cite[Theorem 3.2]{Zqmv} for more details.

Similarly, the analytic continuation of multiple $q$-zeta
functions $\zq(s_1,\dots,s_d)$ can be obtained by using
Euler-Maclaurin summation formula. The major difference between
ordinary MZVs and M$q$ZVs is the appearance of the shifting
operators $\sif_j$ ($1\le j\le d$) defined as follows:
\begin{equation*}
\sif_j\zq(s_1,\dots,s_d):=\zq(s_1,\dots,s_d)+(1-q)\zq(s_1,\dots,s_j-1,\dots,s_d).
\end{equation*}
In general we may iterate the operator and get
\begin{equation}\label{equ:shift}
\sif_j^n\zq(s_1,\dots,s_d)=\sum_{r=0}^n {n\choose r} (1-q)^r
\zq(s_1,\dots,s_j-r,\dots,s_d).
\end{equation}
Using these operators we proved the  analytic continuation of
multiple $q$-zeta functions in \cite{Zqmv}.

The $q$-analogue of MZVs will be called multiple $q$-zeta values
(M$q$ZVs). In this paper we will consider the renormalization
problem for M$q$ZVs motivated by the ideas of Guo and Zhang in
\cite{GuZ}. From physics point of view these values can be
regarded as the quantumization of MZVs. Furthermore these M$q$ZVs
also have number theoretical significance. For instance, it is
well known that $\zeta(0)=-1/2$ and $\zeta(1-2n)=-B_{2n}/(2n)$ for
positive integers $n$ where $B_{2n}$ are Bernoulli numbers defined
by $x/(e^x-1)=\sum_{i=0}^\infty B_ix^i/i!$. What are the right
$q$-analogue of these numbers? It turns out that one of the ways
to find the answer is to consider Rieman $q$ZVs at negative
integers (see \cite[(6)]{KKW}), which shows that the odd indexed
$q$-analogues of Bernoulli numbers are actually nonzero. Is it
possible to generalize Bernoulli numbers to multi-Bernoulli
numbers and their $q$-analogues? Maybe this problem can be solved
when we carry out further studies of the renormalization of MZVs
at negative integers.

The major behavioral difference between MZVs and M$q$ZVs is the
appearance of the shifting operator in the $q$-analogues defined
by \eqref{equ:shift}. As we mentioned in the above it is very
fruitful to study the stuffle relations between MZVs. The
$q$-analogue of this is a little more complicated because of the
shifting operator which can still be handled by setting things up
carefully. The main result of this paper is that we can define the
renormalization of M$q$ZVs  when $s_j$'s are all positive integers
or all non-positive integers such that (i) they coincide with the
M$q$ZV if it is defined originally, (ii) they satisfy a shifted
version of $q$-stuffle relations, and (iii) they become the
renormalization of MZVs defined by Guo and Zhang in \cite{GuZ}
when $\qup$.

To conclude this introduction we remark that currently there are
two ways to order the variables in MVZs and our
$\zeta(s_1,\dots,s_d)$ in this paper is denoted by
$\zeta(s_d,\dots,s_1)$ in \cite{Zqmv}. We change our notation
system because it is more convenient now for the readers to
compare results in this paper to their classical counterparts in
\cite{GuZ} which serves as the major motivation for us.

I would like to thank Li Guo and Bin Zhang for their interest in
this work and many valuable comments of the first draft of this
paper. Thanks are also due to the referees for their many detailed
suggestions which have greatly improved the exposition of this
paper.

\section{The Rota-Baxter algebra and the $q$-stuffle product} Let
$\bfk$ be a subring of $\C$ which is usually taken to be $\R$ or
$\C$. For any fixed $\lambda \in \bfk$ a Rota-Baxter
$\bfk$-algebra of weight $\lambda$ (previously called a Baxter
algebra) is a pair $(R,P)$ in which $R$ is a $\bfk$-algebra and
$P: R \to R$ is a $\bfk$-linear map, such that
\begin{equation}
 P(x)P(y) = P(xP(y)) + P(P(x)y)+ \lambda P(xy),\ \forall x, y \in R.
\label{equ:Ba}
\end{equation}

In this paper we are going to concentrate on the following two examples,
both of which are contained in \cite{GuZ}.
\begin{eg} \label{eg:1} Let $\vep$ be a complex variable such that $\Re(\vep)<0$.
Let $\C\{\!\{\vep,\vep^{-1}\}$ be the algebra of convergent
Laurent series in a neighborhood of $\vep=0$ with at worst finite
order pole at 0. Write $T=-\ln (-\vep)$ which is transcendental
over $\C\{\!\{\vep,\vep^{-1}\}$. Then we can regard
$R:=\C\{\!\{\vep,\vep^{-1}\} [T]$ as the polynomial algebra with
the variable $T$ and with coefficients in
$\C\{\!\{\vep,\vep^{-1}\}$. Let $P$ be the operator on $R$ which
takes the pole part. Then it's not hard to verify that $(R,P)$ is
a Rota-Baxter $\C$-algebra of weight $-1$.
\end{eg}

\begin{eg} \label{eg:Hopf}
Let $\calh$ be a connected filtered Hopf algebra over $\bfk$ (see
\cite[\S2.1]{GuZ} for the definition) and let $(R,P)$ be a
commutative Rota-Baxter algebra of weight $\lambda$. Define the $\bfk$-algebra
$\calr:=\Hom_\bfk(\calh,R)$ of linear maps from $\calh$ to $R$ with the product
compatible with the coproduct of the Hopf algebra $\calh$. Then the operator
$\calp$ on $\Hom(\calh,R)$ defined by $\calp(\calL)=P\circ \calL$
is a Rota-Baxter operator on $\calr$ of weight $\lambda$. This example will be
used in section 2 to define the regularized M$q$ZVs (see \eqref{equ:call}).
\end{eg}

In the rest of this section we will construct one such Hopf
algebra of Example \ref{eg:Hopf}. For any subset $\calz$ of $\C$
closed under addition and shifting by $-1$ we define the
commutative semigroup
\begin{equation}
\frakN(\calz):= \left\{ \wvec{s}{r}\ \Big|\ n\in \Z_{\ge  0},
(s,r)\in \calz\times \R_{>0}\right\}
\end{equation}
with the binary operation given by $\displaystyle \wvec{s}{r}\cdot
\wvec{s'}{r'}=\wvec{s+s'}{r+r'}.$ We will only have two different
choices for $\calz$ in this paper: $\Z$ or $\Z_{\le 0}$. The
reason to require $\calz$ to be closed under shifting by $-1$ is
because of the effects of shifting operators on M$q$ZVs. To study
other renormalization at other poles in the future we need to set
$\calz=\Z+(2 \pi i/\ln q)\Z$ (see \eqref{fSd'}).

Define the $\C$-bilinear pairing $\langle\ \ ,\ \ \rangle$ on the
$\C$-algebra $\C\, \frakN(\calz)$ by
 \begin{equation}\label{equ:pairing}
 \Big\langle \wvec{s}{r}, \wvec{s'}{r'}\Big\rangle
:=\wvec{s+s'}{r+r'}+(1-q) \wvec{s+s'-1}{r+r'}.
\end{equation}
Recall from \cite[\S3.1]{GuZ} that we can define the algebra:
$$\calh_\calz:=\sum_{n\geq 0} \C\, \frakN(\calz)^n$$
where $\frakN(\calz)^0=\{\bfone\}$ is the multiplicative identity
and $\C\, \frakN(\calz)^n$ is the free $\C$-module with basis
$\frakN(\calz)^n$. Then we may equip the $q$-stuffle product,
which is the $q$-analog of the quasi-shuffle product $\ast$ for
MZVs (see \cite[Thm.~2.1]{Ho} or \cite[Thm.~2.2]{GuZ}), on
$\calh_\calz$ as follows: for $\fraka=(a_1, \dots,a_m)\in
\frakN(\calz)^m$ and $\frakb=(b_1,\dots,b_n)\in \frakN(\calz)^n$
we set $\fraka'=\bfone$ if $m=1$ and $\fraka'=(a_2, \dots,  a_m)$
otherwise. Then we define $1\msh \fraka=\fraka\msh 1=\fraka$ and
recursively
 \begin{equation}\label{equ:stuf}
 \fraka \msh \frakb  =
  (a_1,\fraka'\msh \frakb )
  + (b_1, \fraka \msh \frakb')
  + (\langle a_1,b_1\rangle,  \fraka'\msh\frakb')
\end{equation}
where  $\langle a_1,b_1\rangle$ is given by \eqref{equ:pairing}.
It has a connected filtered Hopf algebra structure over $\C$ when
we define the deconcatenation coproduct suitably. If $\vs \msh
\vs'=\sum_{m,n\ge 0} (1-q)^n \vs_{m,n}$ denotes the top component
of $\displaystyle \wvec{\vs}{\vr} \msh \wvec{\vs'}{\vr'}$, then we
proved in \cite[Theorem 5.1]{Zqmv} that
\begin{equation}\label{equ:zqqstuffle}
    \zq(\vs)\zq(\vs')= \zq(\vs \msh \vs'):=
    \sum_{m,n\ge 0} (1-q)^n \zq(\vs_{m,n}).
\end{equation}

\section{Regularized multiple $q$-zeta values}\label{sec:reg}
Let's recall the classical process of renormalization. For
example, let's consider the divergent series $\sum_{n=1}^\infty
n$, which is the series we would get if we tried to plug $s=-1$
into  $\zeta(s)$ by using definition \eqref{equ:oldzeta}. We may
tamper this series by multiplying a controlling factor on each
term: $\sum_{n=1}^\infty ne^{n\vep}$, for some $\vep$ such that
$\Re(\vep)<0$ so that we get a convergent series. By easy
manipulation (see \cite[(34)]{GuZ})
 $$\sum_{n=1}^\infty  n e^{n\vep}=\frac{2}{\vep^2}
 +\sum_{j=0}^\infty  -\frac{B_{j+2}}{j+2}\frac{\vep^j}{j!}.$$
We then call this the \emph{regularized zeta value} at $-1$. To
recover the finite value $\zeta(-1)$ we only need to drop the pole
part $2/\vep^2$ and then take $\vep=0$. This process is called the
\emph{renormalization}. Because there are more than one variable
in multiple zeta functions it turns out that we need to introduce
a concept called ``directional vector'' (see
Definition~\ref{defn:directed}) in the regularization process in
order to get well-behaved regularized values so that the
normalization works as desired.

Turning to our M$q$ZVs, as in the previous section we let $\calz$
be a subset of $\C$ which is closed under addition and shifting by
$-1$ (which will be either $\Z$ or $\Z_{\le 0}$). For $s\in\calz$,
$r>0$, $\Re(\vep)<0$, and $x\in \R$ we first define
$$f(s,r;\vep,x):= \frac{q^{x(s-1)}\exp(r\vep[x]/q^x)}{[x]^s}.$$
Note that the controlling factor becomes $e^{\vep x}$ when $r=1$
and $\qup$. For every vector $\vs=(s_1,\cdots,s_d)\in\calz^d$ and
$\vr=(r_1,\cdots,r_d) \in (\R_{>0})^d$  we now set
\begin{equation}
\Zq\Big(\wvec{\vs}{\vr};\vep,x\Big):=\sum_{n_1>\cdots>n_d>0}
\prod_{j=1}^d f(s_j,r_j;\vep,n_j+x)
\end{equation}
It is clear that $\displaystyle
\Zq\Big(\wvec{\vs}{\vr};\vep,x\Big)$ is also given by the
recursive definition for $\vs=(s_1,\cdots,s_d)\in \calz^d$ and
$\vr=(r_1,\cdots,r_d) \in (\R_{>0})^d$ in \eqref{equ:call2}.
Following the setup of Example \ref{eg:Hopf} we may define the
$\C$-linear map
\begin{align}\label{equ:call}
\calL: \calh_\Z&\lra
R_1:=\left\langle\Zq\Big(\wvec{\vs}{\vr};\vep,x\Big)\left|
  \wvec{\vs}{\vr}\in \frakN\right.\right\rangle    \\
\Big(\wvec{s_1}{r_1},\dots,\wvec{s_d}{r_d}\Big)&\lms
\Zq\Big(\wvec{\vs}{\vr};\vep,x\Big):=Q\Big( f(s_d,r_d;\vep,x)
 \Zq\Big(\wvec{s_1,\cdots,s_{d-1}}{r_1,\cdots,r_{d-1}};\vep,x\Big)\Big),
 \label{equ:call2}
\end{align}
where $Q$ is the summation operator (denoted by $P$ in \cite{EG2,Zud})
$$Q(f)(x)=\sum_{n\geq 1} f(x+n). $$
\begin{defn}\label{defn:regMqZV}
For $\vs\in \calz^d$ and $\vr\in(\R_{>0})^d$ by setting $x=0$ in
$\displaystyle \Zq\Big(\wvec{\vs}{\vr};\vep, x\Big)$ we define
\begin{equation}\label{equ:epsqzeta}
\Zq\Big(\wvec{\vs}{\vr};\vep\Big):=\sum_{k_1>\cdots>k_d>0}
\prod_{j=1}^d \frac{q^{k_j(s_j-1)} \exp(\vep
r_j[k_j]/q^{k_j})}{[k_j]^{s_j}}.
\end{equation}
These values are called the {\em regularized
multiple $q$-zeta values (at $\calz$)}.
\end{defn}
Because of the assumption $\Re(\vep)<0$ we see that
 $\Zq\Big({\displaystyle\wvec{\vs}{\vr}};\vep\Big)$ is well-defined for all $\vs$ and
$\vr$. In particular we don't need to restrict to a MZV-algebra as
constructed in \cite[\S3.2]{EG2} by Ebrahimi-Fard and Guo.
Moreover, when $\qup$ we recover the definition of regularized
general MZV defined in \cite{GuZ}.

\subsection{Regularized $q$-Riemann zeta values}
In this subsection we deal with the $q$ analogue of the Riemann
zeta function first. Taking $d=1$ in Definition \ref{defn:regMqZV}
we find that
$\Zq\Big({\displaystyle\wvec{s}{r}};\vep\Big)=\Zq\Big({\displaystyle\wvec{s}{1}};r\vep)$
so that it suffices to consider
$\displaystyle\Zq(s;\vep):=\Zq\Big(\wvec{s}{1};\vep\Big).$ We will
first put $\calz=\Z$. To study these regularized values we set,
similar to \cite{KKW},
\begin{equation*}
F(x,s;\vep)=\frac{q^{x(s-1)}\exp(\vep q^{-x}[x])}{(1-q)^s [x]^s}=
\frac{q^{x(s-1)}\exp(\vep(q^{-x}-1)/(1-q))} {(1-q^x)^s}, \quad
\Re(\vep)<0 .
\end{equation*}
Then taking derivatives of $F$ with respect to $x$ we get
\begin{align}
F'&(x,s;\vep)= (\ln q)q^{x(s-1)}
\frac{s-1+q^x}{(1-q^x)^{s+1}}\exp(\vep q^{-x}[x])
 -\frac{\vep}{1-q} (\ln q) q^{x(s-2)} \frac{\exp(\vep q^{-x}[x])}{(1-q^x)^s},\notag\\
F''&(x,s;\vep)=(\ln q)^2q^{x(s-1)}
\frac{s(s+1)-3s(1-q^x)+(1-q^x)^2}{(1-q^x)^{s+2}}\exp(\vep q^{-x}[x])\notag\\
-& \frac{\vep}{1-q}(\ln q)^2 q^{x(s-2)}
\frac{2s-3+3q^x}{(1-q^x)^{s+1}}\exp(\vep q^{-x}[x]) +
\Big(\frac{\vep}{1-q}\Big)^2 (\ln q)^2 q^{x(s-3)} \frac{\exp(\vep
q^{-x}[x])}{(1-q^x)^s}. \label{equ:F''}
\end{align}
Let $t=\vep(q^{-x}-1)/(1-q)$, then $dt=-\vep (\ln q)q^{-x}/(1-q)\,
dx$. We get
 \begin{equation} \label{equ:Fintpart}
 \int_{k}^\infty F(x,s;\vep)\, dx  = \frac{1}{\ln q}
 \Big(\frac{\vep}{1-q}\Big)^{s-1}\int_{-\infty}^{\vep(q^{-k}-1)/(1-q)}
  \frac{e^{t} }{t^s} \, dt .
\end{equation}
 When $s=0$ and $k=1$ we have
\begin{equation}
  \label{equ:s=0Fintpart}
  \int_{1}^\infty F(x,0;\vep)\, dx
  = \frac{q-1}{\ln q} \Big(\frac{-1}{\vep}\Big)+
  \sum_{l=0}^\infty \frac{q-1}{\ln q}\cdot\frac{q^{-l-1}}{ -l-1 } \cdot \frac{\vep^l}{l!}.
\end{equation}
Simple computations yields
\begin{align*}
F(1,0;\vep)=&   q^{-1}\exp(\vep/q)= \sum_{l=0}^\infty q^{-l-1}
\frac{\vep^l}{l!}\\
F'(1,0;\vep)=& (-\ln q) q^{-1}+  \sum_{l=1}^\infty
    \frac{\ln q}{1-q} q^{-1-l}(-l-1+q) \frac{\vep^l}{l!},\\
F''(x,0;\vep)=&(\ln q)^2q^{-x}\sum_{l=2}^{\infty}
\Big((q^{-x}[x])^l
 -\frac{3lq^{-x} }{q-1}(q^{-x}[x])^{l-1} + \frac{l(l-1)q^{-2x} }{(q-1)^2}  (q^{-x}[x])^{l-2}
  \Big) \frac{\vep^l}{l!},\\
 \ &\hskip3cm +(\ln q)^2q^{-x}+(\ln q)^2q^{-2x} \frac{3+1-q^x }{1-q}\vep
\end{align*}
Note that for all $0<q<1$ and $\Re(\vep)<0$, $\Zq(0;\vep)=
\sum_{n=1}^\infty F(n,0;\vep)$ converges. By Euler-Maclaurin
summation formula and the analytic continuation of $\zq(s)$ given
by \cite[(12)]{KKW}
\begin{align}
\Zq(0;\vep)
 =&   \int_1^\infty F(x,0;\vep)\, dx +\frac{1}{2}F(1,0;\vep)-
\frac{1}{12} F'(1,0;\vep)- \frac{1}{2}\int_1^\infty\tB_2(x)F''(x,0;\vep)\,dx \nonumber  \\
=&\frac{q-1}{\ln q}\Big(\frac{-1}{\vep}\Big)+\sum_{l=0}^\infty
 \zq(-l)\frac{\vep^l}{l!}, \label{equ:negreg}
\end{align}
where $\tB_2(x)$ is the periodic Bernoulli polynomial defined by
\eqref{equ:bernoulli}. Further, we have (\cite[Ch.~IX,
Misc.~Ex.~12]{WW}) for $k\ge 2$
\begin{equation} \label{equ:perber}
\tB_k(x)=-k!\sum_{n\in \Z\setminus\{0\}} \frac{e^{2\pi inx}}{(2\pi
in)^k},\qquad \tB'_{k+1}(x)=(k+1)\tB_{k}(x).
\end{equation}

To determine the regularized normalization for $\Zq(s;\vep)$ for
positive $s$ we begin with the case $s=1$.  First we have
$$ \int_{1}^\infty F(x,1;\vep)\, dx
 =\frac{-1}{\ln q} \int_{-\frac{\vep}{q}}^\infty \frac{e^{-t} }{t} \, dt.$$
Therefore
\begin{equation}\label{equ:1toinf}
 (1-q) \int_1^\infty F(x,1;\vep)\, dx= \frac{q-1}{\ln q}
 \int_{-\frac{\vep}{q}}^\infty \frac{e^{-t} }{t} \,dt.
\end{equation}
Now for any real number $a>0$ integration by parts yields
$$
\ \int_a^\infty \frac{e^{-t} }{t} \,dt=\Big[e^{-t}\ln
t\Big]_a^\infty+\int_a^\infty   e^{-t} \ln t \, dt
=-\gam- \ln a-\sum_{l=1}^\infty \frac{(-a)^l}{l!l}
$$
where $\gam\approx 0.577216$ is the Euler's $\gam$ constant. Here
we have used the fact that (see \cite{AS} or \cite{Ar})
$$\Gamma'(1)= \int_0^\infty  e^{-t} \ln t \, dt=-\gam.$$
Hence
$$
  (1-q) \int_1^\infty F(x,1;\vep)\, dx=
  \frac{1-q}{\ln q} \left[ \ln \Big(\frac{-\vep}{q}\Big)+\gam +
 \sum_{l=1}^\infty q^{-l}\frac{\vep^l}{l!l} \right]
  $$
Further,
\begin{align*}
F'(1,1;\vep)=& (\ln q)   \frac{q}{(1-q)^2} \exp(\vep/q)
 - \frac{\vep}{1-q} (\ln q) q^{-1} \frac{\exp(\vep/q)}{1-q},\\
F''(x,1;\vep)=&(\ln q)^2
\frac{2-3(1-q^x)+(1-q^x)^2}{(1-q^x)^{3}}\exp(\vep q^{-x}[x])\\
-& \frac{\vep}{1-q}(\ln q)^2 q^{-x}
\frac{-1+3q^x}{(1-q^x)^{2}}\exp(\vep q^{-x}[x]) +
\Big(\frac{\vep}{1-q}\Big)^2 (\ln q)^2 q^{-2x} \frac{\exp(\vep
q^{-x}[x])}{1-q^x}.
\end{align*}
Consequently
\begin{align}
\Zq(1;\vep)= &
 (1-q) \sum_{n=1}^\infty F(n,1;\vep)  \nonumber \\
=&(1-q)\Big( \int_1^\infty F(x,1;\vep)\, dx
 +\frac{1}{2}F(1,1;\vep)-\frac{1}{12} F'(1,1;\vep)
 -\frac{1}{2}\int_1^\infty\tB_2(x)F''(x,1;\vep)\,dx\Big)  \nonumber \\
=& \frac{1-q}{\ln q} \ln (-\vep)+  M(q)+\frac{1-q}{\ln q}
\gam+O(\vep ), \label{equ:posreg}
\end{align}
where we have set
\begin{equation}\label{equ:Mq}
M(q)=q-\frac{1}{2}
 +\frac{q}{12}\frac{\ln q}{q-1}
 - \frac{(1-q)(\ln q)^2}{2}\int_1^\infty\tB_2(x)
 \frac{2-3(1-q^x)+(1-q^x)^2}{(1-q^x)^{3}}\,dx.
\end{equation}
As a comparison we now take a look at the behavior of $\zq(s)$ near
$s=1$. It is clear that
\begin{equation*}
F(x,s,0)= \frac{q^{x(s-1)}}{(1-q^x)^s}.
\end{equation*}
Thus from formula \cite[(12)]{Zqmv} we have, near $s=1$
\begin{align}
\zq(s) =&(1-q)^s\Big(\sum_{n=1}^\infty F(n,s,0)\Big) \nonumber \\
=&(1-q)^s\Big( \int_{1}^\infty F(x,s,0)\, dx
 +\frac{1}{2}F(1,s,0)-\frac{1}{12}F'(1,s,0)
 -\frac{1}{2}\int_1^\infty\tB_2(x)F''(x,s,0)\,dx\Big) \nonumber \\
=&\frac{(q-1)}{(s-1)\ln q}+M(q)+\frac{1-q}{\ln q} \gam+O(s-1).
\end{align}
Taking $\qup$ we have near $s=1$
\begin{equation}\label{equ:zetaat0}
 \zeta(s)=\frac{1}{s-1}+\lim_\qup M(q)-\gam+O(s-1).
\end{equation}
On the other hand, by applications of Euler-Maclaurin summation
formula (\cite[p.~531]{Knopp} or \cite[7.21]{WW}), for all
integers $k$, $l\ge 1$, we get
$$\sum_{j=1}^k \frac{1}{j}=\ln k
 +\gam+\frac{1}{2k}-\sum_{r=1}^l\frac{B_{2r}}{2r}
 k^{-2r}+\int_k^\infty \frac{\tB_{2l+1}(x)}{x^{2l+2}}\,dx.$$
Since $B_2=1/6$, putting $k=l=1$, using equation
\eqref{equ:perber} and integration by parts once in \eqref{equ:Mq}
we get $\lim_\qup M(q)=\gam$ which is consistent with
\eqref{equ:zetaat0}.

Now we can prove the following result:
\begin{thm} \label{thm:Riemanncase}
The value of $\Zq\Big({\displaystyle \wvec{s}{r}};\vep\Big)$ at
$s=1-n$ ($n\in\Z_{>0}$) is
\begin{equation}\label{equ:negr}
 \Zq(1-n;r\vep)=\frac{q-1}{\ln q}\left(\frac{-1}{r\vep}\right)^n n!
 +\sum_{l=0}^\infty \zq(1-n-l)\frac{(r\vep)^l}{l!}.
\end{equation}
The value of $\Zq\Big({\displaystyle \wvec{s}{r}};\vep\Big)$ at
$s=n\in\Z_{>0}$ is
\begin{equation}\label{equ:posr}
 \Zq(n;r\vep)=\frac{q-1}{\ln q}
 u_n(r\vep)+\Big(M(q)+\frac{1-q}{\ln q} \gam\Big)\frac{(r\vep)^{n-1}}{(n-1)!}
 +\sum_{l=0,l\ne n-1}^\infty \zq(n-l)\frac{(r\vep)^l}{l!}.
\end{equation}
where $u_{n+1}(\vep)=\frac{ \vep^n}{n!}(H_n-\ln (-\vep))$ for
$n\ge 0$, $H_0=0$ and $H_n=1+\frac12+\cdots+\frac1n$ for $n\ge 1.$
\end{thm}
\begin{proof}  Equation \eqref{equ:negr} follows from
\eqref{equ:negreg} by taking derivatives because
$\frac{d}{d\vep}\Zq(s;\vep)=\Zq(s-1;\vep)$. When $s=1$ equation
\eqref{equ:posr} follows from \eqref{equ:negreg} and
\eqref{equ:posreg} by integration and Abel's Theorem because
$\frac{d}{d\vep}\Zq(1;\vep)=\Zq(0;\vep).$ The rest follows
immediately.
\end{proof}

\subsection{The range of regularized M$q$ZV}
We now turn to the general M$q$ZVs. Although the proof of
Theorem~\ref{thm:zlaurent} below is similar to the proof of
\cite[Theorem 3.3]{GuZ} some new phenomena arise because of the
shifting principle for M$q$ZV. Let
 $$\C\{\!\{\vep,\vep^{-1}\}=\left\{\sum_{n=N}^\infty
  a_n\vep^n\Bigg| a_n\in \C, N\in \Z\right\}$$
be the algebra of Laurent series, regarded as a subalgebra of the
algebra of (the germs of) meromorphic functions in a neighborhood
of $\vep=0$ with at worst finite order poles at $0$. Choose $\ln
\vep$ to be analytic on $\C\backslash (-\infty,0]$. Observe that
the analytic function $\ln (-\vep)$ on $\C\backslash [0,\infty)$
is transcendental over $\C\{\!\{\vep,\vep^{-1}\}$ by \cite[Lemma
3.1]{GuZ} and hence there is a natural algebra injection (see
\cite[(19)]{GuZ})
  \begin{equation}\label{equ:u}
u:\C\{\!\{\vep,\vep^{-1}\} [-\ln(-\vep)]\lra
 \C[T]\{\!\{\vep,\vep^{-1}\}
\end{equation}
sending $-\ln(-\vep)$ to $T$. This provides an identification of
$\C\{\!\{\vep,\vep^{-1}\} [-\ln(-\vep)]$ as a subalgebra of
$\C[T]\{\!\{\vep,\vep^{-1}\}$.

\begin{thm}\label{thm:zlaurent}
\emph{(a)} For any $\vs\in (\Z_{\le 0})^d$ and $\vr\in
(\R_{>0})^d$,
 $\Zq\Big({\displaystyle \wvec{\vs}{\vr}};\vep\Big) \in\C\{\!\{\vep,\vep^{-1}\}$.

\emph{(b)} For any  $\vs\in \Z^d$ and $\vr\in (\R_{>0})^d$,
$\Zq\Big({\displaystyle \wvec{\vs}{\vr}};\vep\Big)\in
\C\{\!\{\vep,\vep^{-1}\}[-\ln (-\vep)]$.
\end{thm}
\begin{proof}
(a)  The key is the following computation of the tail of the
$\Zq(s;\vep)$ for $s\in\Z_{\le 0}$. First define
$$\xi_i(\vep):=\sum_{k> i}\frac{\exp(\vep[k]/q^k)}{q^k}
 =\sum_{k>0}\frac{\exp(\vep[k+i]/q^{k+i})}{q^{k+i}}.$$
Since $[k+i]=[k]+q^k[i]$ we get
$$\xi_i(\vep)=\frac{\exp(\vep[i]/q^i)}{q^{i}} \cdot
 \sum_{k>0}\frac{\exp(\vep[k]/q^{k+i})}{q^{k}}
 =\frac{\exp(\vep[i]/q^i)}{q^{i}} \Zq(0;\vep/q^i) .$$
Therefore  for $s\in \Z_{\le 0}$ and $t=-s$,
\begin{align}
 \ & \xi_i(t;\vep):=\sum_{k> i}
 \frac{q^{k(-t-1)}\exp(\vep[k]/q^k)}{[k]^{-t}}=
 \Big(\frac{d}{d\vep}\Big)^t \xi_i(\vep)
 =  \sum_{j=0}^t {t\choose j} \frac{[i]^j\exp(\vep[i]/q^i)}{q^{i(j+1)}}
 \Big(\frac{d}{d\vep}\Big)^{t-j} \Zq(0;\vep/q^i) \nonumber\\
 =&\frac{q-1}{\ln q}\sum_{j=0}^t \frac{t!}{j!}\Big(\frac{ -1 }{\vep}\Big)^{1+t-j}
  \frac{[i]^j\exp(\vep[i]/q^i)}{q^{ij}}
 +\sum_{j=0}^t {t\choose j} \sum_{l=t-j}^\infty
   \frac{[i]^j\exp(\vep[i]/q^i)}{q^{i(j+1+l)}}
 \frac{\zq(-l)\vep^{l-t+j}}{(l-t+j)!} .  \label{equ:xi}
\end{align}

Now we first prove (a) by induction on the length $d$. The case
$d=1$ corresponds to the Riemann $q$-zeta function which has been
dealt with in the last section. Suppose (a) is true for length
$d-1$ ($d>1$) and let $\vs=(s_1,\dots,s_d).$ Let $\vs_{\hat{i}}$
denote $\vs$ with its $i$-th component removed. Let $\bfe_i$
denote the $i$-th unit vector of length $d-1$ with 1 at the $i$-th
component. By definition \eqref{equ:epsqzeta} and the shifting
operator \eqref{equ:shift} for every $i$ such that $1\le i\le d$,
\begin{align} \notag
 \ & \Zq\Big(\wvec{\vs}{\vr};\vep\Big) =\sum_{k_1>\cdots>k_d>0}
   \prod_{m=1}^{d}
\frac{q^{k_m(s_m-1)} \exp(\vep r_m[k_m]/q^{k_m})}
 {[k_m]^{s_m}}  \\
 =&\sum_{k_1>\cdots>k_{i-1}} \sum_{k_{i+1}>\cdots>k_d>0}
  \prod_{m=1,m\ne i}^{d}
\frac{q^{k_m(s_m-1)} \exp(\vep r_m[k_m]/q^{k_m})}
 {[k_m]^{s_m}} \notag  \\
\ &\hskip2cm \cdot \Big[ \xi_{k_{i+1}}(-s_i;r_i\vep) -
\xi_{k_{i-1}}(-s_i; r_i\vep)-
  \frac{q^{k_{i-1}(s_i-1)}
 \exp(\vep  r_i[k_{i-1}]/q^{k_{i-1}})}{[k_{i-1}]^{s_i}}\Big] \notag\\
=& - \sif_{i-1}
\Zq\Big(\wvec{\vs_{\hat{i}}+s_i\bfe_{i-1}}{\vr_{\hat{i}}+r_i
\bfe_{i-1}};\vep\Big)+\frac{q-1}{\ln
q}\sum_{j=0}^{-s_i}\frac{(-s_i)!}{j!} \Big(\frac{-1
}{r_i\vep}\Big)^{1-s_i-j}
 \Big[ \Zq\Big(\wvec{\vs_{\hat{i}}-j\bfe_i}{\vr_{\hat{i}}+r_i \bfe_i};\vep\Big)
  - \Zq\Big(\wvec{\vs_{\hat{i}}-j\bfe_{i-1}}{\vr_{\hat{i}}+r_i \bfe_{i-1}};\vep\Big)\Big] \notag\\
\ &+\sum_{j=0}^{-s_i}  {-s_i\choose j}\sum_{l=-s_i-j}^\infty
 \Big[\sif_i^{l+1}\Zq\Big(\wvec{\vs_{\hat{i}}-j\bfe_i}{\vr_{\hat{i}}+r_i \bfe_i};\vep\Big)
  - \sif_{i-1}^{l+1}\Zq\Big(\wvec{\vs_{\hat{i}}-j\bfe_{i-1}}{\vr_{\hat{i}}+r_i
  \bfe_{i-1}};\vep\Big)\Big]
  \frac{\zq(-l) (r_i\vep)^{l+s_i+j}}{(l+s_i+j)!}.
\label{equ:denom}
\end{align}
Here if $i=1$ then the terms with $\bfe_{i-1}$ are $0$. If $i=d$
then the terms with $\bfe_i$ are $0$ and the first term becomes
$\displaystyle\Zq\Big(\wvec{s_d}{r_d};\vep\Big)
\Zq\Big(\wvec{\vs_{\hat{d}}}{\vr_{\hat{d}}};\vep\Big)$.
By induction assumption we see that $\Zq\Big({\displaystyle
\wvec{\vs}{\vr}};\vep\Big)\in\C\{\!\{\vep,\vep^{-1}\}$.

Now we prove (b) by induction on length $d$ again.

\noindent {\bf Case 1.} $s_i\leq 0$ for some $1\leq i\leq d$. Then
the proof of (a) above carries over word for word here.

\noindent {\bf Case 2.} Suppose $s_i>0$ for all $1\le i\le d$. We
use induction on the sum $s:=\sum_{i=1}^{d} s_i$. Clearly $s\ge
d$. If $s=d$, then $s_i=1$ for $1\le i\le d$. Thus
$$\Zq'\Big(\wvec{\vs}{\vr}; \vep\Big)
 =\sum r_i \Zq\Big(\wvec{\vs-\bfe_i}{\vr}; \vep\Big) \in\C\{\!\{\vep,\vep^{-1}\}$$
by (a). Integrating we get $\Zq\Big({\displaystyle
\wvec{\vs}{\vr}};\vep\Big)
\in\C\{\!\{\vep,\vep^{-1}\}[-\ln(-\vep)]$. The general case
follows from the fact that $\C\{\!\{\vep,\vep^{-1}\}[-\ln(-\vep)]$
is closed under integration.
\end{proof}

Later on, for several times we are going to need the special case
of \eqref{equ:denom} when $i=1$ and $i=d$. For convenience we list
them as
\begin{cor} \label{cor:special1d}
Let $d\ge 2$ be a positive integer.
For any $\vs\in (\Z_{\le 0})^d$ and $\vr\in
(\R_{>0})^d$,
\begin{equation}\label{equ:denom1}
  \aligned
  \Zq\Big(\wvec{\vs}{\vr};\vep\Big)  =&  \frac{q-1}{\ln
q}\sum_{j=0}^{-s_1}\frac{(-s_1)!}{j!} \Big(\frac{-1
}{r_1\vep}\Big)^{1-s_1-j}
  \Zq\Big(\wvec{s_2-j,s_3\dots,s_d}{r_1+r_2,r_3,\dots,r_d};\vep\Big)
  \\
\ &+\sum_{j=0}^{-s_1}  {-s_1\choose j}\sum_{l=-s_1-j}^\infty
 \sif_1^{l+1}\Zq\Big(\wvec{s_2-j,s_3\dots,s_d}{r_1+r_2,r_3,\dots,r_d};\vep\Big)
  \frac{\zq(-l) (r_1\vep)^{l+s_1+j}}{(l+s_1+j)!},
\endaligned
\end{equation}
and
\begin{equation}\label{equ:denomd}
\aligned
  \Zq\Big(\wvec{\vs}{\vr};\vep\Big)  =&\Zq\Big(\wvec{s_d}{r_d};\vep\Big)
\Zq\Big(\wvec{s_1,\dots,s_{d-1}}{r_1,\dots,r_{d-1}};\vep\Big)\\
 &- \frac{q-1}{\ln q}\sum_{j=0}^{-s_d}\frac{(-s_d)!}{j!}
\Big(\frac{-1 }{r_d\vep}\Big)^{1-s_d-j}
\Zq\Big(\wvec{s_1,\dots,s_{d-2},s_{d-1}-j}{r_1,\dots,r_{d-2},r_{d-1}+r_d};\vep\Big) \\
\ &-\sum_{j=0}^{-s_d}  {-s_d\choose j}\sum_{l=-s_d-j}^\infty
\sif_{d-1}^{l+1}\Zq\Big(\wvec{s_1,\dots,s_{d-2},s_{d-1}-j}{r_1,\dots,r_{d-2},r_{d-1}+r_d};\vep\Big)
  \frac{\zq(-l) (r_d\vep)^{l+s_d+j}}{(l+s_d+j)!}.
\endaligned
\end{equation}
\end{cor}

The next result tells us some very useful information about the
general shape of the coefficients of the Laurent series
$\displaystyle\Zq\Big(\wvec{\vs}{\vr};\vep\Big)$. This will be
used crucially in the proof of the existence of the
renormalizations of M$q$ZVs at non-positive integers.
\begin{cor} \label{cor:denomr}
Let $\vs\in (\Z_{\le 0})^d$ and $\vr\in (\R_{>0})^d$. Then the
coefficient $c_i$ of the Lauren series of
$\displaystyle\Zq\Big(\wvec{\vs}{\vr};\vep\Big)=\sum c_i\vep^i$ is
an $\R$-linear combination of rational function of the form
$P(\vr)/Q(\vr)$ in $r_1,\dots,r_d$, where $P$, $Q$ have no common
factors. Moreover, both $P$ and $Q$ are homogeneous polynomials in
$r_1,\dots,r_d$, $Q$ is the product of linear factor
$r_1+\cdots+r_d$ (and $r_1$ if $d=2$) with repetition allowed, and
$\deg(P/Q):=\deg(P)-\deg(Q)=i$.
\end{cor}
\begin{proof} When $d=1$ this follow from
Theorem~\ref{thm:Riemanncase}. Suppose the claim in the corollary
is true when the length of the vector $\vs$ is $<d$ for some $d\ge
2$. It follows from \eqref{equ:denom1} and induction assumption
that the only possible factors in the denominator of $P/Q$ are of
the form $(r_1+\cdots+r_d)^j$ and $r_1^k$ for some $j,k\in \Z_{\ge
0}$. But by \eqref{equ:denomd} we see that the only possibilities
are $(r_1+\cdots+r_d)^j$, $(r_1+\cdots+r_{d-1})^k$, and $r_d^l$.
Hence if $d>2$ then none of $r_1$, $r_1+\cdots+r_{d-1}$ or $r_d$
can appear as a factor in the denominator of $c_i$. But if $d=2$
then both $r_1+r_2$ and $r_1$ could appear (and indeed they do by
the following formula \eqref{equ:uZq2}.)

The statements about the degrees are clear from \eqref{equ:denom}
by induction assumption.
\end{proof}

\section{Renormalization of M$q$ZV}
Theorem~\ref{thm:zlaurent} together with the map $u$ of
\eqref{equ:u} shows that there is an algebra homomorphism
$$\aligned
\uZq: \calh_{\Z}\lra  &\C[T]\{\!\{\vep,\vep^{-1}\}\\
\wvec{\vs}{\vr}\longmapsto &
u\Big(\Zq\Big(\wvec{\vs}{\vr};\vep\Big)\Big)
\endaligned
$$
which restricts to an algebra homomorphism
 $$\uZq: \calh_{\Z_{\le 0}}\to\C\{\!\{\vep,\vep^{-1}\}.$$
Note that only when there is a positive $s_i$ can $\uZq$ really
differ from $\Zq$. It is well-known that the Birkhoff
decomposition (see \cite[Theorem II.5.1]{Ma}) yields two maps
${\tilde Z}_{q-}$ and ${\tilde Z}_{q+}$ such that $\uZq={\tilde
Z}_{q-}^{-1} \star {\tilde Z}_{q+}.$ The properties of this
decomposition implies the following result immediately.

\begin{prop}\label{prop:renormz}
The map ${\tilde Z}_{q+}: \calh_{\Z}\to \C[T][\![\vep]\!]$ is an
algebra homomorphism which restricts to an algebra homomorphism
${\tilde Z}_{q+}: \calh_{\Z_{\le 0}}\to \C[\![\vep]\!]$.
\end{prop}

To write down ${\tilde Z}_{q+}$ explicitly we need the following
definition which is slightly different from
\cite[Definition~3.7]{GuZ}.
\begin{defn} \label{de:part}
Let $\Pi_d$ be the set of increasing sequences $i_0=0<
i_1<\cdots<i_p=d$. For $1\leq j\leq p$, define the \emph{partition
vectors} of $\vs\in \C^d$ from the sequence $(i_1,\cdots,i_p)$ to
be the vectors $\vs^{(j)}:=(s_{i_{j-1}+1},\cdots,s_{i_{j}})$,
$1\leq j\leq p$.
\end{defn}
The following explicit formula for the renormalization of
regularized M$q$ZV is the $q$-analog of \cite[Theorem~3.8]{GuZ}
which follows from \cite[Theorem II.5.1]{Ma}, which, in turn, is
built upon the idea of \cite[Theorem 4]{CK}.
\begin{prop}\label{prop:Z+}
Let $P$ be the operator sending a Laurent series to its pole part:
$P(\sum_{n\geq -N}^\infty a_n \vep^n)=\sum_{n=-N}^{-1}a_n \vep^n$
and let $\check{P}=-P$. Then for any $\vs\in \Z^d$ and $\vr\in
(\R_{>0})^d$
\begin{equation}
\label{equ:Z+1} {\tilde Z}_{q+}\Big(\wvec{\vs}{\vr}\Big)
=\sum_{(i_1,\cdots,i_p)\in\Pi_d}\hspace{-.6cm} \tilde{P}\llp
\uZq\lp\wvec{\vs^{(p)}}{\vr^{(p)}}\rp
    \cdots \check{P}\llp \uZq\lp\wvec{\vs^{(2)}}{\vr ^{(2)}}\rp
    \check{P}\llp \uZq\lp\wvec{\vs ^{(1)}}{\vr ^{(1)}}\rp\rrp\,
    \rrp \cdots \rrp.
\end{equation}
\end{prop}

We are going to use the map ${\tilde Z}_{q+}$ to define the
renormalization of M$q$ZVs. But before doing so we recall that at
the beginning of section \ref{sec:reg} we mention that the
multiple variable cases are different from the single variable
case, just like the situation where a function of two variables
can have all the directional derivatives at some point but yet is
not differentiable there. Such phenomenon won't happen to
functions of one variable. The renormalization process is
essentially a limit process just like taking the derivatives. The
behavior of multiple zeta functions at poles are not so bad in
that if we take appropriate ``paths'' to renormalize then we can
produce values compatible with both the stuffle relations and the
original values if they are defined originally. All the above
remarks are still valid for the $q$-analogues. The appropriate
``paths'' in our case is given by
Definition~\ref{defn:noshiftMqZV} and
Definition~\ref{defn:shifted} later in this section. We first
define the directional version of these as follows:
\begin{defn} \label{defn:directed}
 For $\vs\in\Z^d$ and $\vr\in (\R_{>0})^d$, the {\em
renormalized directional M$q$ZV} is defined by
\begin{equation*}
    \zq\Big(\wvec{\vs}{\vr}\Big)
    = \lim_{\vep\to 0} {\tilde Z}_{q+}\Big(\wvec{\vs}{\vr} \Big),
\end{equation*}
and $\vr$ is called the {\em directional vector}.
\end{defn}

\begin{cor} \label{cor:stuffle}
The renormalized directional M$q$ZVs satisfy the $q$-stuffle
relations
\begin{equation}\label{equ:zqdirqstuffle}
    \zq\Big(\wvec{\vs}{\vr}\Big)\zq\Big(\wvec{\vs'}{\vr'}\Big)
    = \zq\Big(\wvec{\vs}{\vr}\msh \wvec{\vs'}{\vr'}\Big):=
    \sum_{n\ge 0} (1-q)^n \sum_{m\ge 0}\zq\Big(\wvec{\vs_{m,n}}{\vr_{m,n}}\Big)
\end{equation}
where ${\displaystyle \wvec{\vs}{\vr}}\msh {\displaystyle
\wvec{\vs'}{\vr'}}=
  \sum_{n\ge 0, m\ge 0} (1-q)^n {\displaystyle \wvec{\vs_{m,n}}{\vr_{m,n}}}$.
\end{cor}
\begin{proof}
It follows directly from Definition~\ref{defn:directed} by
Proposition~\ref{prop:renormz}.
\end{proof}

The following proposition is the $q$-analogue of
\cite[Corollary~4.13]{GuZ}. This and the length three case provide
us a hint at the general shapes of the renormalization of M$q$ZVs
at non-positive integers which will be given in
Proposition~\ref{prop:mn}

\begin {prop} \label{prop:z2formula}
For $s_1$, $s_2 \le 0$, $r_1$, $r_2>0$, set $t=1-s_1-s_2$. Then
\begin{align*}
\zq\Big(\wvec{s_1, s_2}{r_1,r_2}\Big)=& \frac{q-1}{\ln q}\left[
\frac{\zq(s_1+s_2-1)}{s_1-1}+\sum_{j=0}^{-s_1} {-s_1\choose j}
\sum_{l=t}^\infty
 {l+1\choose l-t}(1-q)^{l-t}
  \Big(\frac{-r_1}{r_1+r_2}\Big)^{l+s_1+j}
  \frac{\zq(-l)}{l+s_1+j}\right]\\
\ &+\sum_{j=0}^{-s_1} {-s_1\choose j}
 \sum_{i=0}^{1-j-s_1}{1-j-s_1\choose i}(1-q)^i \zq(s_2-j-i)
   \zq(j+s_1-1).
\end{align*}
\end{prop}
\begin{proof}
 From equation \eqref{equ:denom1} in Corollary \ref{cor:special1d} we
 have
\begin{align*}
&\uZq\Big(\wvec{s_1,s_2}{r_1,r_2}\Big)=\frac{q-1}{\ln q}
\sum_{j=0}^{-s_1}\frac{(-s_1)!}{j!}\Big(\frac{-1}{r_1\vep}\Big)^{1-s_1-j}
 \Zq\Big(\wvec{s_2-j}{r_1+r_2};\vep\Big) \\
\ &+\sum_{j=0}^{-s_1} {-s_1\choose j} \sum_{l=-s_1-j}^\infty
 \sum_{i=0}^{l+1}{l+1\choose i}(1-q)^i\Zq\Big(\wvec{s_2-j-i}{r_1+r_2};\vep\Big)
  \frac{\zq(-l) (r_1\vep)^{l+s_1+j}}{(l+s_1+j)!}
\end{align*}
By Theorem \ref{thm:Riemanncase} we get
\begin{align}
\uZq &\Big(\wvec{s_1,s_2}{r_1,r_2}\Big)=\left(\frac{q-1}{\ln
q}\right)^2
\sum_{j=0}^{-s_1}\frac{(-s_1)!}{j!}\Big(\frac{-1}{r_1\vep}\Big)^{1-s_1-j}
 \left(\frac{-1}{(r_1+r_2)\vep}\right)^{j-s_2+1}(j-s_2)! \notag\\
 \ &+\frac{q-1}{\ln q}\sum_{j=0}^{-s_1}\frac{(-s_1)!}{j!}\Big(\frac{-1}{r_1\vep}\Big)^{1-s_1-j}
 \sum_{k=0}^\infty \zq(s_2-j-k)\frac{(
 (r_1+r_2)\vep)^{k}}{(k)!}\label{equ:uZq2}\\
\ &+\frac{q-1}{\ln q}\sum_{j=0}^{-s_1} {-s_1\choose j}
\sum_{l=-s_1-j}^\infty
 \sum_{i=0}^{l+1} {l+1\choose i}(1-q)^i \frac{(r_1\vep)^{l+s_1+j}}{(-(r_1+r_2)\vep )^{j+i-s_2+1}}
 \frac{(j+i-s_2)!\zq(-l) }{(l+s_1+j)!} \notag\\
\ &+\sum_{j=0}^{-s_1} {-s_1\choose j} \sum_{l=-s_1-j}^\infty
 \sum_{i=0}^{l+1}{l+1\choose i} \sum_{k=0}^\infty \zq(s_2-j-i-k)\frac{(
 (r_1+r_2)\vep)^{k}}{(k)!} (1-q)^i
  \frac{\zq(-l) (r_1\vep)^{l+s_1+j}}{(l+s_1+j)!}\notag
\end{align}
By Definition \ref{defn:directed} and Proposition~\ref{prop:Z+} we
have
 \begin{align*}
\zq\Big(\wvec{s_1, s_2}{r_1,r_2}\Big) =\lim_{\vep\to
0}\left[{\tilde P}\Big(\Zq\Big(\wvec{s_2}{r_2};\vep\Big) {\check
P}\Big(\Zq\Big(\wvec{s_1}{r_1};\vep\Big) \Big) \Big) + {\tilde
P}\Big(\Zq\Big(\wvec{s_1,s_2}{r_1,r_2};\vep\Big)\Big)\right].
\end{align*}
It follows from Theorem \ref{thm:Riemanncase} and the above
expression for
$\displaystyle\Zq\Big(\wvec{s_1,s_2}{r_1,r_2};\vep\Big)$ that
 \begin{align*}
&\zq\Big(\wvec{s_1, s_2}{r_1,r_2}\Big)=\frac{q-1}{\ln q}
\left(-\Big(\frac{-r_2}{r_1}\Big)^{1-s_1}\frac{\zq(s_1+s_2-1)}{1-s_1}+
\sum_{j=0}^{-s_1}\frac{(-s_1)!}{j!}\Big(\frac{r_1+r_2}{-r_1}\Big)^{1-s_1-j}
\frac{\zq(s_1+s_2-1)}{(1-s_1-j)!}\right)\\
\ &\quad +\frac{q-1}{\ln q}\sum_{j=0}^{-s_1} {-s_1\choose j}
\sum_{l=1-s_1-s_2}^\infty
 {l+1\choose l+s_1+s_2-1}(1-q)^{l+s_1+s_2-1}
  \Big(\frac{-r_1}{r_1+r_2}\Big)^{l+s_1+j}
  \frac{\zq(-l)}{l+s_1+j}\\
\ &\quad +\sum_{j=0}^{-s_1} {-s_1\choose j}
 \sum_{i=0}^{1-j-s_1}{1-j-s_1\choose i}(1-q)^i \zq(s_2-j-i)
   \zq(j+s_1-1).
\end{align*}
A simple combinatorial formula quickly reduces this to the formula
in the proposition.
\end{proof}

We are now ready to define the renormalization of M$q$ZVs.
\begin{defn}  \label{defn:noshiftMqZV}
For $\vs \in (\Z_{>0})^d\cup (\Z_{\le 0})^d$ define
\begin{equation}\label{eq:gqmzv}
\gzq(\vs)= \lim_{\delta \to 0^+}
\zq\Big(\wvec{\vs}{|\vs|+\delta}\Big),
\end{equation}
where, for $\vs=(s_1,\cdots,s_d)$ and $\delta\in \R_{>0}$, we
write $|\vs|=(|s_1|,\cdots,|s_d|)$ and
$|\vs|+\delta=(|s_1|+\delta,\cdots,|s_d|+\delta).$ These values
$\gzq(\vs)$ are called the {\em renormalized M$q$ZV} of the
$d$-tuple $q$-zeta function $\zq(u_1,\cdots,u_d)$ at $\vs$.
\end{defn}
In order to deal with the $q$-stuffle relations we need the
shifted version of the about definition.
\begin{defn} \label{defn:shifted}
Let $\vf=(f_1,\dots,f_l)$ be a \emph{binary} vector with entries
equal to either $0$ or $1$. Let $\vs\in(\Z_{\le 0})^l$ or
$\vs\in(\Z_{>0})^l$ such that $|\vs|-\vf\in\Z_{\ge 0}^l$ (i.e., if
$s_i=0$ then $f_i=0$). We call such a binary vector $\vf$ a
shifting vector of $\vs$. We define the \emph{shifted
renormalization of $\zq(\vs)$ by $\vf$} as the limit
\begin{equation}\label{equ:shifted} \gzq^{
\vf}(\vs):=\lim_{\delta\to 0^+}
 \zq\Big(\wvec{\vs}{\delta+|-\vs+(\delta-1)\vf|}\Big).
\end{equation}
For example, $\displaystyle\gzq^{(1,0,0)}(-1,-2,0)=\lim_{\delta\to
0^+} \zq\Big(\wvec{-1,-2,0}{2\delta,2+\delta,\delta}\Big)$, where
$2\delta$ appears because of the stuffings at the first position
which is hinted by the shifting vector $(1,0,0)$.
\end{defn}

We will compute the limit \eqref{eq:gqmzv} in
Theorem~\ref{thm:sallpos} and show that the limits
\eqref{equ:shifted} exist for $\vs \in (\Z_{>0})^d$ in
Theorem~\ref{thm:shiftedspos}, for $\vs \in(\Z_{\le 0})^d$ in
Theorem~\ref{thm:gzeta}. First we consider the case of Riemann
$q$-zeta function.
\begin{thm} Let $\gam$ be Euler's gamma constant and define $M(q)$ by \eqref{equ:Mq}. Then
 $$\gzq(1)=\frac{q-1}{\ln q}T + M(q)+\frac{1-q}{\ln q}\gam$$
and, for integers $s>1$, $\gzq(s)$ is the
usual Riemann $q$-zeta value $\zq(s)$ defined by the series \eqref{equ:qzeta}.
If $s=-l$ is a non-positive integer then
\begin{equation}
  \label{eq:zqneg}
  \gzq(-l)=\gz_q(-l)=(1-q)^{-l}\left\{
\sum_{r=0}^l(-1)^{r}\binom{l}{r}\frac1{q^{l+1-r}-1}
+\frac{(-1)^{l+1}}{(l+1)\ln q}\right\}.
\end{equation}
\end{thm}
\begin{proof}
The expression of $\gz_q(-l)$ is given by \cite[(60)]{KKW}. The
rest follows from Theorem \ref{thm:Riemanncase}.
\end{proof}

Our primary goals are to show that our definition of
renormalizations of M$q$ZVs are well-define, these values agree
with the usual M$q$ZVs whenever the usual values are defined, they
satisfy the $q$-stuffle relations, and they become
renormalizations of MZVs when $\qup$.

\begin{thm}\label{thm:sd>1}
Let $\vs=(s_1,\cdots,s_d)\in\Z^d_{>0}$ such that $s_1>1$. Let
$\vr\in (\Z_{>0})^d$ be an arbitrary vector. Then
$\zq\Big({\displaystyle \wvec{\vs}{\vr}}\Big)=\zq(\vs)$ which is
independent of the choice of $\vr$. Here $\zq(\vs)$ denotes the
usual definition of multiple $q$-zeta values. In particular,
$\gzq(\vs)=\zq(\vs)$ satisfy the $q$-stuffle relation.
\end{thm}
\begin{proof} By definition $\Zq\Big({\displaystyle \wvec{\vs}{\vr}};\vep\Big)$
converges uniformly for $\vep\in(-\infty,0]$ and therefore
continuous as a function of $\vep$. In particular it is a power
series and $\lim_{\vep\to 0}\Zq\Big({\displaystyle
\wvec{\vs}{\vr}};\vep\Big)=\Zq\Big({\displaystyle
\wvec{\vs}{\vr}};0\Big)=\zq(\vs)$.
\end{proof}

Now we consider the divergent case $s_1=1$. For
$f(\vep),g(\vep)\in \C[T][\![\vep]\!]$, denote $f(\vep)=g(\vep)+
O(\vep)$ if $g(\vep)-f(\vep)\in \vep \C[T][\![\vep]\!]$.

\begin {lem} \label{lem:X}
For $c>0$ set $X=\frac {q-1}{\ln q}(\ln c+T)$. Let $\vs \in
(\Z_{>0})^l$ where either $l=0$ or $s_1>1$. Let $\vr \in
(\R_{>0})^l$ and $\vone_d =(1,1,\cdots, 1)\in \Z^d$. Then
\begin{equation}\label{equ:111}
\uZq\Big(\wvec{\vone_d,\vs}{c\vone_d,\vr}\Big)= P_{d,\vs}(X)+
O(\vep),
\end{equation}
where $P_{d,\vs}(X)$ is a degree $d$ polynomial in $X$ with
leading coefficient $\zq(\vs)/d!$ (by convention, if $l=0$ we set
$\zq(\vs)=1$). Moreover, $P_{d,\vs}(X)$ is independent of $\vr$.
\end{lem}

\begin {proof} We prove the lemma by induction similar to
that of \cite[Lemma 4.4]{GuZ}. Notice that the terms produced by
shifting will produce polynomials of degrees less than that of the
leading term.

When $d=0$ the lemma follows from the proof of
Theorem~\ref{thm:sd>1} as $\vep \C[\![\vep]\!]\subset\vep
\C[T][\![\vep]\!]$. We now fix $d=1$ and prove the lemma by
induction on $l$. When $l=0$ the lemma readily follows from
\eqref{equ:posr}. Assume equation \eqref{equ:111} is true when the
length of the vector $\vs$ is $\le l$ for $l\ge 0$. Then by
\eqref{equ:posr} and Theorem~\ref{thm:sd>1} we have
\begin{multline*}
  (X+O(\vep))\uZq\Big(\wvec{\vs}{\vr}\Big)-\uZq\Big(\wvec{1,\vs}{c,\vr}\Big)
  =u\Big(\Zq\Big(\wvec{1}{c};\vep\Big)\cdot\Zq\Big(\wvec{\vs}{\vr};\vep\Big)
  -\Zq\Big(\wvec{1,\vs}{c,\vr}\vep\Big)\Big)\\
= \sum_{j=1}^{l-1}
\uZq\Big(\wvec{s_1,\cdots,s_j,1,s_{j+1},\cdots,s_l}{r_1,\cdots,r_j,c,r_{j+1},\cdots,r_l}\Big)
+\sum _{j=1}^l \uZq\Big(\wvec{\vs+ \bfe_j}{\vr+
c\bfe_j}\Big)+(1-q)\sum _{j=1}^l \uZq\Big(\wvec{\vs}{\vr+
c\bfe_j}\Big)=O(\vep)
\end{multline*}
where $\bfe_j$ is the $j$-th unit vector of length $l$. Hence by
induction assumption,
 $$\uZq\Big(\wvec{1,\vs}{c,\vr}\Big)=\zq(\vs)X+O(\vep).$$

Assume equation \eqref{equ:111} is true for every $l>0$ when the
length of the vector $(1,\dots,1)$ is $\le d$ for $d\ge 1$. Then
by \eqref{equ:posr}
 $$ \uZq\Big(\wvec{\vone_d,\vs}{c\vone_d, \vr}\Big)\cdot
 \uZq\Big(\wvec{1}{c}\Big) =(P_d(X)+ O(\vep))(X+O(\vep))
 =f_{d+1}(X)+ O(\vep)$$
for some polynomial $f_{d+1}(X)$ of degree $d+1$ with leading
coefficient $\zq (\vs)/d!$, independent of $\vr$. On the other
hand, by the $q$-stuffle relation
 \begin{multline*}
 \uZq\Big(\wvec{\vone_d,\vs}{c\vone_d,\vr}\Big)\cdot\uZq\Big(\wvec{1}{c}\Big)
=(d+1) \uZq\Big(\wvec{\vone_{d+1},\vs}{c\vone_{d+1},\vr}\Big)\\
 + \sum_{j=1}^{l-1} \uZq\Big(\wvec{\vone_d,s_1,\cdots,s_j,1,s_{j+1},\cdots,s_l}
 {c\vone_d,r_1,\cdots,r_j,c,r_{j+1},\cdots,r_l}\Big)
 +\sum _{j=1}^l \uZq\Big(\wvec{\vone_d, \vs+\bfe_j}{c\vone_d,\vr+c\bfe_j}\Big)
 +(1-q)\sum _{j=1}^l \uZq\Big(\wvec{\vone_d,\vs}{c\vone_d,\vr+ c\bfe_j}\Big).
\end{multline*}
By induction assumption
 $$\uZq\Big(\wvec{\vone_{d+1},\vs}{c\vone_{d+1},\vr}\Big)=P_{d+1,\vs}(X)+
 O(\vep)$$
where $P_{d+1,\vs}(X)$ is some polynomial of degree $d+1$ with
leading coefficient $\zq (\vs)/(d+1)!$, independent of $\vr$. This
completes the proof of the lemma.
\end{proof}
Now we can show
\begin {thm}\label{thm:sallpos}
For $\vs \in (\Z_{>0})^m$,  $\gzq(\vs)=\zq\Big({\displaystyle
\wvec{\vs}{\vs}}\Big)$.\label{thm:gregmzv}
\end{thm}
\begin{proof} In view of Theorem \ref{thm:sd>1} we only need to
consider the case $s_1=1$. Assume then $\vs=(\vone_d,\vs')$ where
$\vs'=(s_1,\dots,s_l)$ and $s_1>1$. Observe that for any substring
$\vs''$ of $\vs$ and substring $\vr''$ of $\vr$ by the above lemma
$\Zq\Big({\displaystyle \wvec{\vs''}{\vr''}}\Big)$ has no pole
part. This means that in equation \eqref{equ:Z+1} of
Proposition~\ref{prop:Z+} there is only one non-trivial term which
gives the finite part
 $$\Zq_+\Big({\displaystyle \wvec{\vs}{\vr}}\Big)=
 (Id-P)\Big(\Zq\Big({\displaystyle \wvec{\vs}{\vr}},\vep\Big)\Big)=P_{d,\vs'}(X)+O(\vep).$$
Hence
 $$\zq\Big({\displaystyle\wvec{\vs}{\vr}}\Big)=P_{d,\vs'}(X).$$
By Definition~\ref{defn:noshiftMqZV}
 $$\gzq(\vs)=\lim_{c\to 1,\vr'\to \vs'}
 \zq\Big({\displaystyle\wvec{\vone_d,\vs'}{c\vone_d,
 \vr}}\Big)=P_{d,\vs'}\Big(\frac{(q-1)T}{\ln q}\Big)=
 \zq\Big({\displaystyle\wvec{\vs}{\vs}}\Big).$$
\end{proof}
However, when $s_1=1$ the $q$-stuffle relation are not exactly
preserved under renormalization because of shifting. For example,
Proposition~\ref{prop:renormz} implies that ${\tilde Z}_{q+}$ is
an algebra homomorphism $\calh_{\Z}\lra \C[T][\![\vep]\!]$ with
$${\tilde Z}_{q+}\Big({\wvec{s_1}{s_1}}, \cdots, {\wvec{s_k}{s_k}}\Big)\Big|_{\vep=0}
 =\zeta\lp\wvec{\vs}{\vs}\rp=\gzq(\vs).$$
Thus for $c>0$ we have
$$\aligned
\uZq\lp\wvec{1}{c}\rp^2=\uZq\lp\wvec{1}{c} \msh
\wvec{1}{c}\rp=&2\uZq\lp\wvec{1,1}{c,c}\rp+
\uZq\lp\wvec{2}{2c}\rp+(1-q)\uZq\lp\wvec{1}{2c}\rp\\
=&2\uZq\lp\wvec{1,1}{c,c}\rp+
\uZq\lp\wvec{2}{2c}\rp+(1-q)\uZq\lp\wvec{1}{c}\rp+\frac{(1-q)^2\ln
2}{\ln q}
\endaligned$$
by Theorem \ref{thm:Riemanncase}. This implies
$\gzq^{(1)}(1)=\gzq(1)+(1-q)\ln 2/\ln q$ and
 $$\aligned
 \gzq(1\msh 1)=&2\gzq(1,1)+\gzq(2)+(1-q)\gzq(1)\\
 \gzq(1)^2=&2\gzq(1,1)+\gzq(2)+(1-q)\gzq^{(1)}(1).\endaligned$$
Of course, when $\qup$ these two values are both equal to
$\bar{\zeta}(1)^2=T^2$.
\begin {thm}\label{thm:shiftedspos}
Let $\vs \in (\Z_{>0})^l$ with $s_1>1$, then for any positive
integer $d$ the limit of $\gzq^{\vf}(\vone_d,\vs)$ in
\eqref{equ:shifted} exists. Moreover all the shifted
renormalization of M$q$ZVs at all positive arguments satisfy the
$q$-stuffle relations in the following sense. Let $\vs_1 \in
(\Z_{>0})^k, \vs_2 \in (\Z_{>0})^\ell$ and assume there are $W_t$
vectors of length $k+\ell-t$ produced by the $q$-stuffle  with $t$
stuffings and let $\vf_j$ be the binary vector representing the
stuffing positions of the $j$-th such vector:
$$ \vs_1\msh \vs_2=\sum_{t=0}^{\min(k,\ell)}
 (1-q)^t\sum_{j=1}^{W_t} (\vec{v}_j-\vf_j).$$
Then
$$\gzq(\vs_1)\gzq(\vs_2)=\sum_{t=0}^{\min{k,\ell}}
 (1-q)^t\sum_{j=1}^{W_t} \gzq^{\vf_j}(\vec{v}_j-\vf_j).$$
\end{thm}
\begin{proof} By the very definition of the shifted
renormalization and the $q$-stuffle relation it suffices to show
the first part of the theorem, namely the existence of the limit
in  \eqref{equ:shifted}. Indeed, if the leading component in $\vs$
is greater than 1 then $\gzq$ is independent of the directional
vector by Theorem~\ref{thm:sd>1}.

In the following we consider the shifted renormalization
$\gzq^{\vf}(\vone_d,\vs)$ only. Let the vector $\vr=
(r_1,\dots,r_{k+l})=(\vone_d,\vs)+\vf$ be the shifted directional
vector of length $l+d$. Then
 \begin{equation}\label{equ:indregval}
\Zq\lp\wvec{\vone_d,\vs}{\vr};\vep\rp =\sum_{k_1>\cdots>k_d}
\prod_{j=1}^d\frac{\exp(r_j\vep[k_j]/q^{k_j})}{[k_j]} \sum_{k_d>n_1>\cdots>n_l}\\
\prod_{j=1}^l\frac{q^{n_j(s_j-1)}\exp(r_j\vep[n_j]/q^{n_j})}{[n_j]^{s_j}}.
\end{equation}
We want to show this series is \emph{good} (used just in this
proof) in the sense that it is in $\C[\![\vep]\!][\ln(-\vep)]$ and
is finite when $r_1,\cdots,r_d$ are equal to either 1 or 2. We
prove this by induction on $d$. If $d=1$ this follows from
\eqref{equ:111} in Lemma \ref{lem:X}. Assume the series in
\eqref{equ:indregval} is good when the length of $(1,\dots,1)$ in
front of $\vs$ is less than $d$ for some $d\ge 2$. Set as before
 $$F(x,s;\vep)=\frac{q^{x(s-1)}\exp(\vep q^{-x}[x])}{(1-q)^s [x]^s}=
 \frac{q^{x(s-1)}\exp(\vep(q^{-x}-1)/(1-q))} {(1-q^x)^s}. $$
As in the proof of Theorem \ref{thm:Riemanncase} we can use
Euler-Maclaurin summation formula to find a closed expression for
\begin{multline}\label{equ:Fs}
 \sum_{k>i}F(k,1;\vep)=\sum_{k>i}\frac{\exp(\vep[k]/q^{k})}{(1-q)[k]}\\
=\int_1^\infty F(x+i,1;\vep)\, dx
 +\frac{1}{2}F(i+1,1;\vep)-\frac{1}{12} F'(i+1,1;\vep)
 -\frac{1}{2}\int_1^\infty\tB_2(x)F''(x+i,1;\vep)\,dx.
\end{multline}
The following integral
$$ \int_{1}^\infty F(x+i,1;\vep)\, dx
 =\frac{-1}{\ln q} \int_{-\frac{\vep[i+1]}{q^{i+1}}}^\infty \frac{e^{-t} }{t} \, dt.$$
can be evaluated using \eqref{equ:1toinf} by the substitution
$\vep\to \vep[i+1]/q^i$. Therefore we get
\begin{multline}\label{equ:integral}
(1-q)\int_{1}^\infty F(x+i,1;\vep)\, dx=  \frac{1-q}{\ln q} \left[
\ln \Big(\frac{-\vep[i+1]}{q^{i+1}}\Big)+\gam +
 \sum_{l=1}^\infty \frac{\vep^l[i+1]^l}{l!l q^{l(i+1)}} \right]\\
 =  \frac{1-q}{\ln q} \left[ \ln(-\vep)-\ln q
+\ln\Big(\frac{[i]}{q^i}\Big)+\ln\Big(1+\frac{q^i}{[i]}[1]\Big)+\gam
+
 \sum_{l=1}^\infty \frac{\vep^l}{l!l q^l}\cdot
  \Big(\frac{[i]}{q^i}+[1]\Big)^l\right].
\end{multline}
Let $i=k_2$ in the above two formulas \eqref{equ:Fs} and
\eqref{equ:integral}. It is straightforward to see that the middle
two terms of \eqref{equ:Fs} both contribute to good series by
Theorem~\ref{thm:sd>1} and the induction assumption. The last term
of  \eqref{equ:Fs} can be handled similarly using \eqref{equ:F''}
after we notice that $\tB_2(x)$ is bounded by $1/6$ from an easy
computation from the series expansion \eqref{equ:perber} and the
fact that $\zeta(2)=\pi^2/6$. So the main divergence term of
\eqref{equ:indregval} when $\vep$ is near zero comes from the
first term of \eqref{equ:Fs}. Thus by \eqref{equ:integral} it is
enough to show that all the sums below are good:
\begin{align}\label{equ:ind1}
\ & \sum_{k_2>\cdots>k_d}
\ln\Big(\frac{[k_2]}{q^{k_2}}\Big)\prod_{j=2}^d\frac{\exp(r_j\vep[k_j]/q^{k_j})}{[k_j]}
\sum_{k_d>n_1>\cdots>n_l} \cdots,\\
\ & \sum_{k_2>\cdots>k_d}
 \frac{q^{mk_2}}{[k_2]^m}\prod_{j=2}^d\frac{\exp(r_j\vep[k_j]/q^{k_j})}{[k_j]}
\sum_{k_d>n_1>\cdots>n_l} \cdots, m\in \Z_{\ge0 } \label{equ:ind2}\\
\ & \sum_{k_2>\cdots>k_d}
\frac{\vep^l[k_2]^m}{q^{mk_2}}\prod_{j=2}^d\frac{\exp(r_j\vep[k_j]/q^{k_j})}{[k_j]}
\sum_{k_d>n_1>\cdots>n_l} \cdots, 1\le m\le l.\label{equ:ind3}
\end{align}
In \eqref{equ:ind2} if $m=0$ then this follows from induction and
if $m>0$ then it follows from Theorem~\ref{thm:sd>1}. For
\eqref{equ:ind3}, we can mimic the argument for finding
$\Zq(s,\vep)$ at non-positive integers by differentiating
\eqref{equ:integral} and see immediately that \eqref{equ:ind3} is
actually in $\C[\![\vep]\!]$ and $r_2$ can appear in the
denominator only in the form of some pure power. So the series in
\eqref{equ:ind3} is good by induction.

The most difficult one is \eqref{equ:ind1} in which case we again
apply Euler-Maclaurin summation formula using the following
modified version of $F(x,s;\vep)$:
\begin{equation*}
G_j(x,s;\vep)=\ln^j\Big(\frac{[x]}{q^{x}}\Big)\frac{q^{x(s-1)}\exp(\vep
q^{-x}[x])}{(1-q)^s [x]^s}.
\end{equation*}
Then \eqref{equ:integral} changes to
\begin{align}\notag
\int_{1}^\infty G_j(x+i,1;\vep)\, dx=&
 \int_{-\frac{\vep[i+1]}{q^i}}^\infty (\ln t)^j t^{-1} e^{-t}  \,dt\\
 =& -\frac{1}{j+1} \ln^{k+1}\Big(-\frac{\vep[i+1]}{q^i}\Big)  \exp\Big(\frac{\vep[i+1]}{q^i}\Big)
 +\frac{1}{j+1}\int_{-\frac{\vep[i+1]}{q^i}}^\infty (\ln t)^{j+1} e^{-t}
 \,dt.\label{equ:integralG}
\end{align}
This integral can be treated similarly as \eqref{equ:integral} by
using polygamma functions (see \cite{AS}) which is closely related
to
$$f_k(s):= \int_0^\infty (\ln t)^k  t^{1-s} e^{-t}  \,dt=\Gamma^{(k)}(s).$$
In particular, to proceed by induction we need the fact that all
the values of $f_k(1)$ are finite which are in fact $\Q$-linear
combinations of the weight $k$ products of $\zeta(n)$'s and the
Euler constant $\gamma$, where we take the weight of $\gamma$ to
be 1. For example, the first few values are $f_1(1)=-\gamma,\
f_2(1)=\zeta(2)+\gamma^2 , \
f_3(1)=-2\zeta(3)-3\gamma\zeta(2)-\gamma^3,$ and
$f_4(1)=\frac{27}{2}\zeta(4)+8\gamma\zeta(3)+6\gamma^2\zeta(2)+\gamma^4.$
This shows that the sum in \eqref{equ:ind1} lies in
$\C[\![\vep]\!]$ and the only form that $r_j$'s can occur in the
denominator of any coefficient of the series is in a factor of
some form $r_{i_1}+\cdots+r_{i_t}$. Note that  $r_j$'s can also
occur in logarithms in the form $\ln(r_j)$ which implies that we
can take  $r_j$'s to be either 1 or 2. This finishes the proof of
the theorem.
\end{proof}

\section{Renormalization of M$q$ZVs at nonpositive integers}
The following result is straightforward:
\begin{thm}\label{thm:lz}
For $\vs\in (\Z_{<0})^d$, we have
$$ \gzq(\vs)=\zq\Big(\wvec{\vs}{-\vs}\Big)
=\lim_{\vr\to -\vs} \zq\Big(\wvec{\vs}{\vr}\Big).$$
\end{thm}
\begin{proof}
This follows from Corollary \ref{cor:denomr},
Definition~\ref{defn:directed} and
Definition~\ref{defn:noshiftMqZV}.
\end{proof}

In the rest of the section we are going to prove
\begin {thm}\label{thm:gzeta}
Let $\vs\in (\Z_{\le 0})^d$ and $\vf$ be a shifting vector of
$\vs$. Then the limit in \eqref{equ:shifted}
$$ \gzq^{
\vf}(\vs):=\lim_{\delta\to 0^+}
 \zq\Big(\wvec{\vs}{\delta+|-\vs+(\delta-1)\vf|}\Big)
$$
exists.
\end {thm}

The major complication of the proof of Theorem~\ref{thm:gzeta}
arises from the possibility that $s_i$ can be 0 or $-1$. To prove
the theorem we need some more information of $\displaystyle{\tilde
Z}_{q+}\Big(\wvec{\vs}{\vr}\Big)$.

\begin{defn} \label{defn:sym}
Let $\vr\in (\R_{>0})^d$. Let $d$ and $j$ be two positive integers
such that $j\le d.$ Let $N(\vr)=r_i+r_{i+1}+\cdot+r_d$,
$D(\vr)=r_j+r_{j+1}+\cdots+r_d$ for $1\le j<i\le d$. Then we say
$N(\vr)^w/D(\vr)^v$ is a \emph{tailored fraction} for any
nonnegative integers $v$ and $w$ such that $w\ge v$. For
convenience, we regard all polynomials of $\vr$ as tailored
fractions too. If a rational function $P/Q$ is a product of
tailored fractions as above then we call it a \emph{tailored
rational function}. Note that $\deg(P/Q)\ge 0$ always holds. If
every coefficient of a Lauren series $\sum c_i\vep^i$ is an
$\R$-linear combination of tailored rational functions then we say
the Lauren series is a {\em tailored Laurent series}.
\end{defn}

\begin{prop} \label{prop:mn}
Suppose $m$ and $n$ are two integers and $d$ is a positive
integer. Let $\vs\in (\Z_{\le 0})^d$ and $\vr\in (\R_{>0})^d$. For
$(i_1,\cdots,i_p)\in \Pi_d$, let $\vec{s}^{(j)}, 1\leq j\leq p$ be
the partition vectors of $(s_2,\dots,s_d)$ from $(i_1,\cdots,i_p)$
in Definition~\ref{de:part} and similarly define $\vec{r}^{(j)}$,
$1\leq j\leq p$. Let
\begin{equation}\label{equ:Zm+}
{\tilde Z}_{q+}^{n,m}\Big(\wvec{\vs}{\vr}\Big)=\sum_{i\ge 0}
c_i\vep^i :=\sum_{(i_1,\cdots,i_p)\in\Pi_d}\hspace{-.6cm}
\tilde{P}\llp \vep^n \uZq\lp\wvec{\vs^{(p)}}{\vr^{(p)}}\rp
    \cdots \check{P}\llp \uZq\lp\wvec{\vs^{(2)}}{\vr ^{(2)}}\rp
    \check{P}\llp  \vep^m\uZq\lp\wvec{\vs ^{(1)}}{\vr ^{(1)}}\rp\rrp\,
    \rrp \cdots \rrp.
\end{equation}
Then for every $i$ the coefficient $c_i$ can be expressed as an
$\R$-linear combination of rational functions $P/Q$ of
$r_1,\cdots,r_d$ with $P$ and $Q$ having no common factors such
that the following must hold:

(a) $P$ and $Q$ are homogeneous polynomials such that
$\deg(P/Q)+n+m=i$.

(b) Either $Q$ is a constant or every factor of $Q$ has the form
$r_j+\cdots+r_d$ for some $j=1,\dots,d$. So $Q$ is uniquely
determined if we require $Q$ to be a monic polynomial with respect
to $r_d$. Then either $Q=1$ or $Q$ is a product of linear factors
of the form $r_j+\cdots+r_d$ ($j=1,\dots,d$) with possible
repetitions.

(c) If $m=0$ then $r_d$ divides $Q$ only if $n>1-s_d$.

(d)  If $m=0$ and $n\le 0$ then $P/Q$ is tailored.
\end{prop}
\begin{proof} (a) By Corollary~\ref{cor:denomr} for $1\le t\le p$
every coefficient $a_i$ of
${\displaystyle\uZq\lp\wvec{\vs^{(t)}}{\vr ^{(t)}}\rp=}\sum_i a_i
\vep^i$ is a homogeneous rational function $P/Q$ of
$r_1,\cdots,r_d$ such that $\deg(P/Q)+i=0$. Now the product of two
such terms
$$ (P_1/Q_1)\vep^i \cdot (P_2/Q_2)\vep^j=(P_3/Q_3)\vep^{i+j}$$
still satisfies $\deg(P_3/Q_3+i+j=0$. Part (a) follows easily.

(b) We proceed by induction on $d$. If $d=1$ then this follows
immediately from Theorem \ref{thm:Riemanncase}. Now assume part
(b) is true when the length of the vector $\vs$ is $d-1$ for some
$d\ge 2.$ Let $\vs\in (\Z_{\le 0})^d$. We have the disjoint union
\begin{align*}
  \Pi_{d}=&\{(1,i_1+1,\cdots,i_p+1)\ \big|\ (i_1,\cdots,i_p)\in \Pi_{d-1}\}  \\
& \cup
    \{(i_1+1,i_2+1,\cdots,i_p+1)\ \big|\ (i_1,\cdots,i_p)\in
    \Pi_{d-1}\}.
\end{align*}
Note that $i_p=d-1$ in the above by definition. In the rest of the
proof of (b), for $(i_1,\cdots,i_p)\in \Pi_{d-1}$ let
$\vec{s}^{(j)}, 1\leq j\leq p$ be the partition vectors of
$(s_2,\dots,s_d)$ from $(i_1+1,\cdots,i_p+1)$ in
Definition~\ref{de:part} and similarly define $\vec{r}^{(j)}$,
$1\leq j\leq p$. Putting
$\sum_{\vi}=\sum_{(i_1,\cdots,i_p)\in\Pi_{d-1}}$ we have
 \begin{align*}
{\tilde Z}_{q+}^{n,m}(\wvec{\vs}{\vr})
 =& \sum_{\vi} \tilde{P}\llp
\vep^n \uZq\lp\wvec{\vs^{(p)}}{\vr^{(p)}}\rp
    \cdots \check{P}\llp \uZq\lp\wvec{\vs^{(1)}}{\vr ^{(1)}}\rp
    \check{P}\llp  \vep^m\uZq\lp\wvec{s_1}{r_1}\rp\rrp\,
    \rrp \cdots \rrp\\
&+\sum_{\vi} \tilde{P}\llp \vep^n
\uZq\lp\wvec{\vs^{(p)}}{\vr^{(p)}}\rp
    \cdots \check{P}\llp \uZq\lp\wvec{\vs^{(2)}}{\vr ^{(2)}}\rp
    \check{P}\llp  \vep^m \uZq\lp\wvec{s_1,\vs^{(1)}}{r_1,\vr ^{(1)}}\rp\rrp\,
    \rrp \cdots \rrp.
\end{align*}
By Theorem \ref{thm:Riemanncase} we have
$$ \check{P}\llp \vep^m \uZq\lp\wvec{s_1}{r_1}\rp\rrp=
\begin{cases} 0 \quad &\text{if }m>-s_1, \\
 {\displaystyle \frac{q-1}{\ln q} \frac{A(r_1)(-s_1)!}{\vep^{1-(m+s_1)}}  -\sum_{i=0}^{-m-1}
 \zq(s_1-i) \frac{r_1^i\vep^{m+k}}{i!} }\quad &\text{if }m\le -s_1,
 \end{cases}$$
where the sum is vacuous if $m\ge 0$ and
\begin{equation}\label{equ:A}
    A(r_1):=-\Big(\frac{-1}{r_1}\Big)^{1-s_1}.
\end{equation}
Write $k=i_1+1$. Then by \eqref{equ:denom1} we get
 \begin{multline*}
\vep^m \uZq\lp\wvec{s_1,\vs^{(1)}}{r_1,\vr ^{(1)}}\rp
 =\frac{q-1}{\ln q} \frac{B(r_1)(-s_1)!}{\vep^{1-(m+s_1)}} \\
 +\sum_{j=0}^{-s_1}  {-s_1\choose j}\sum_{l=-s_1-j}^{\infty}
 \sif_1^{l+1}\Zq\Big(\wvec{s_2-j,s_3,\dots,s_k}{r_1+r_2,r_3,\dots,r_k};\vep\Big)
  \frac{\zq(-l) (r_1\vep)^{l+s_1+j}}{(l+s_1+j)!\vep^{-m}},
\end{multline*}
where
\begin{equation}\label{equ:B}
    B(r_1):=\sum_{j=0}^{-s_1}\frac{1}{j!}
\Big(\frac{-1 }{r_1}\Big)^{1-s_1-j}
 \Zq\Big(\wvec{s_2-j,s_3,\dots,s_k}{r_1+r_2,r_3,\dots,r_k};\vep\Big)\vep^j.
\end{equation}
By induction assumption, the linear factors of $Q$ can only be
$r_1$, or $r_j+\cdots+r_d$ ($j=1,2,\dots,d$). Let us exclude the
possibility of $r_1$. It suffices to prove that
$$\lim_{r_1\to 0} \left(A(r_1)\uZq\lp\wvec{\vs^{(1)}}
 {\vr ^{(1)}}\rp+B(r_1)\right)\  \text{ is finite.}$$
Set $\displaystyle
f(x)=\Zq\Big(\wvec{s_2,s_3,\dots,s_k}{x+r_2,r_3,\dots,r_k};\vep\Big)$.
Put $D=d/dx$ and $t=-s_1\ge 0$. By L'Hopital's rule
 \begin{multline*}
\lim_{x\to 0} f(0)A(x)+B(x)
 =(-1)^{t}\lim_{x\to 0} \frac{\displaystyle f(0)-
 \sum_{j=0}^{t}\frac{1}{j!} (-x)^j
  D^j f(x)}{x^{t+1}}\\
 = (-1)^t\lim_{x\to 0} \frac{\displaystyle
 \sum_{j=0}^{t}\frac{1}{j!}
 (-x)^j D^{j+1}f(x)-\sum_{j=1}^{t} \frac{1}{(j-1)!}  (-x)^{j-1}
 D^j f(x)}{(t+1)x^t}
 = \frac{(D^{t+1}f)(0)}{(t+1)!} <\infty
\end{multline*}
after cancellations in the numerator.

\begin{rem}\label{rem:rd}
By exactly the same argument we see that $r_d$ can not appear in
the denominator of \eqref{equ:denomd} when $d\ge 2$.
\end{rem}

(c) and (d). Assume $m=0$. We use induction on $d$ again. Theorem
\ref{thm:Riemanncase} yields the case $d=1$ easily for both (c)
and (d). Suppose part (d) holds when $\vs$ has length $d-1$ for
some $d\ge 2.$ Now let $\vs\in (\Z_{\le 0})^d$. Observe that we
have another disjoint union of $\Pi_{d}$ given by
\begin{align*}
  \Pi_{d}=&\{(i_1,\cdots,i_p,d)\ \big|\ (i_1,\cdots,i_p)\in \Pi_{d-1}\}  \\
& \cup
    \{(i_1,\cdots,i_{p-1},i_p+1)\ \big|\ (i_1,\cdots,i_p)\in
    \Pi_{d-1}\}.
\end{align*}
Similar to part (b), in the rest of the proof, for any positive
integer $t$ and $(i_1,\cdots,i_p)\in \Pi_t$, let $\vs^{(j)}, 1\leq
j\leq p$ be the partition vectors of $(s_1,\dots,s_t)$ from
$(i_1,\cdots,i_p)$ in Definition~\ref{de:part}. Define $\vr^{(j)}$
in the same fashion. Put
$\sum_{\vi(t)}=\sum_{(i_1,\cdots,i_p)\in\Pi_t}$. Then
 \begin{align*}
{\tilde Z}_{q+}^{n,0}(\wvec{\vs}{\vr})
 =& \sum_{\vi(d-1)} \tilde{P}\llp
\vep^n \uZq\lp\wvec{s_d}{r_d}\rp\check{P}\llp
\uZq\lp\wvec{\vs^{(p)}}{\vr^{(p)}}\rp
    \cdots \check{P}\llp\uZq\lp\wvec{\vs^{(1)}}{\vr^{(1)}}\rp
    \rrp \cdots \rrp\rrp\\
&+\sum_{\vi(d-1)} \tilde{P}\llp \vep^n
\uZq\lp\wvec{\vs^{(p)},s_d}{\vr^{(p)},r_d}\rp
    \check{P}\llp \uZq\lp\wvec{\vs^{(p-1)}}{\vr^{(p-1)}}\rp
    \cdots \check{P}\llp \uZq\lp\wvec{\vs^{(1)}}{\vr ^{(1)}}\rp\rrp\,
     \cdots\rrp \rrp.
\end{align*}
For ease of notation, in the rest of the proof we write
$R_i=r_i+\dots+r_d$, $k=i_{p-1}+1$, and $\displaystyle X=
\check{P}\llp \uZq\lp\wvec{\vs^{(p-1)}}{\vr ^{(p-1)}}\rp \cdots
\check{P}\llp \uZq\lp\wvec{\vs^{(1)}}{\vr ^{(1)}}\rp\rrp \cdots
\rrp$. Using \eqref{equ:denomd} we get
\begin{align}
\ &{\tilde Z}_{q+}^{n,0}(\wvec{\vs}{\vr})
 =\sum_{\vi(d-1)}
  \tilde{P}\llp\vep^n \uZq\lp\wvec{s_d}{r_d}\rp\check{P}\llp
\uZq\lp\wvec{\vs^{(p)}}{\vr^{(p)}}\rp X\rrp+
\vep^n \uZq\lp\wvec{\vs^{(p)},s_d}{\vr^{(p)},r_d}\rp X\rrp \label{equ:induc0}\\
=&\sum_{\vi(d-1)}
 \tilde{P}\llp\vep^n \uZq\lp\wvec{s_d}{r_d}\rp{\tilde P}\llp
\uZq\lp\wvec{\vs^{(p)}}{\vr^{(p)}}\rp X\rrp \rrp\label{equ:induc1}\\
 &- \frac{q-1}{\ln q}\sum_{j=0}^{-s_d} \frac{(-s_d)!}{j!}\sum_{\vi(d-1)}
 \tilde{P}\llp\vep^n \Big(\frac{-1 }{r_d\vep}\Big)^{1-s_d-j}
\Zq\Big(\wvec{s_k,\dots,s_{d-2},s_{d-1}-j}{r_k,\,\dots\,
,r_{d-2},\, \, \ R_{d-1}\, \, };\vep\Big)X \rrp \label{equ:induc2}\\
\ &-\sum_{j=0}^{-s_d}  {-s_d\choose j}\sum_{l=-s_d-j}^\infty
  \sum_{\vi(d-1)} \hskip-.2cm\tilde{P}\llp\vep^n   \frac{\zq(-l) (r_d\vep)^{l+s_d+j}}{(l+s_d+j)!}
  \sif_{d-1}^{l+1}\Zq\Big(\wvec{s_k,\dots,s_{d-2},s_{d-1}-j}{r_k,\,\dots\,
,r_{d-2},\, \, \ R_{d-1}\, \, };\vep\Big) X \rrp.
\label{equ:induc3}
\end{align}
By Theorem \ref{thm:Riemanncase} and Remark~\ref{rem:rd} we see
that $r_d$ can appear essentially in the denominator of some term
in \eqref{equ:induc0} only if $n>1-s_d$ which proves (c). As a
matter of fact, we don't need induction assumption to prove (c) so
we will use (c) freely in what follows.

Suppose now $n\le 0$. By part (b) every denominator in
\eqref{equ:induc1} can only have linear factors of the form $r_d$
or $R_j-r_d$ for some $j=1,\dots,d-1$. By part (b) and (c) we know
that these terms will be cancelled out in the end so
\eqref{equ:induc1} doesn't contribute to any un-tailored terms. In
fact, negative powers of $R_j-r_d$ can not really appear by the
reasoning in the next paragraph.

As the exponent of $\vep$ in front of $\Zq$ is now negative in
\eqref{equ:induc2} every term has $r_d$ in the denominator so that
all of these will be cancelled by corresponding terms in
\eqref{equ:induc2} (note that there is no negative powers of $r_d$
in \eqref{equ:induc3}). In fact by part (b) the only linear
factors of denominators that can appear are of the form $R_j$ in
\eqref{equ:induc2} and of the form either $R_j-r_d$ or $r_d$ in
\eqref{equ:induc1}. So to cancel terms with negative $r_d$-power
no factors $R_j-r_d$ can appear in denominators of
\eqref{equ:induc1} and no factors $R_j$ can appear in denominators
of \eqref{equ:induc2}. (Actually, positive powers $R_j$ and
$R_j-r_d$ must appear so that $r_d=R_j-(R_j-r_d)$ is produced on
the numerator to cancel negative powers of $r_d$). This shows that
the sum of \eqref{equ:induc1} and \eqref{equ:induc2} is tailored.

Let's turn to \eqref{equ:induc3}. If $l+j$ is small, say $l+j\le
-s_d-n$ then by induction assumption the Laurent series is
tailored. When $l+j$ increase gradually it seems that un-tailored
terms may appear. We will show this actually can not happen.
Roughly speaking, in order for some linear factor
$r_j+r_{j+1}+\cdots+r_d$ to appear in the denominator, the power
of $\vep$ in front of $\Zq$ in \eqref{equ:induc3} has to be large
by part (c). But this will produce high powers  of $r_d$ which
results in tailored fractions.

It is easy to see that in \eqref{equ:induc3} the highest power of
$r_{d-1}+r_d$ that can appear in the denominator is $1-s_{d-1}+j$
by part (c). This takes place if $l+s_d+j+n>1-s_{d-1}+j$.
Similarly, we can expand \eqref{equ:induc3} by iteratively using
the formulas from \eqref{equ:induc0} to \eqref{equ:induc3} and see
that the highest power of $R_i$ that can appear in the denominator
is $1-s_i+j_{d-i}$ where $j_{d-i}$ is the corresponding summation
index in \eqref{equ:induc2} and \eqref{equ:induc3}. Indeed
iterating once we see that  for fixed $j$ and $l$ the essential
part of \eqref{equ:induc3} becomes:
\begin{align}
\ & r_d^{l+s_d+j}\sum_{\vi(d-1)}\tilde{P}\llp\vep^{n+l+s_d+j}
\Zq\Big(\wvec{s_k,\dots,s_{d-2},s_{d-1}-j}{r_k,\,\dots\,
,r_{d-2},\, \, \ R_{d-1}\, \, };\vep\Big)
 X \rrp  \label{equ:induc4} \\
 =&\sum_{\vi(d-2)}
 r_d^{l+s_d+j}\tilde{P}\llp\vep^{n+l+s_d+j}  \uZq\lp\wvec{s_{d-1}-j}{R_{d-1}}\rp{\tilde P}\llp
\uZq\lp\wvec{\vs^{(p)}}{\vr^{(p)}}\rp X\rrp \rrp\label{equ:induc11}\\
 & \aligned
 &- \frac{q-1}{\ln q}\sum_{j_2=0}^{j-s_{d-1}} \frac{(j-s_{d-1})!}{j_2!}\sum_{\vi(d-2)}
\\
 \ & \hskip1cm   r_d^{l+s_d+j}\tilde{P}\llp \vep^{n+l+s_d+j}
 \Big(\frac{-1}{R_{d-1}\vep}\Big)^{1+j-s_{d-1}-j_2}
 \Zq\Big(\wvec{s_k,\dots,s_{d-3},s_{d-2}-j_2}{r_k,\,\dots\,
,r_{d-3},\, \ \ R_{d-2}\, \ };\vep\Big)X
 \rrp \endaligned  \label{equ:induc5}\\
 \ &\aligned
 &-\sum_{j_2=0}^{j-s_{d-1}} {j-s_{d-1}\choose j_2}\sum_{l_2=j-s_{d-1}-j_2}^\infty
  \sum_{\vi(d-2)}\frac{\zq(-l_2)
}{(l_2+s_d-j+j_2)!} \cdot\\
\ & \hskip1cm r_d^{l+s_d+j}\tilde{P}\llp\vep^{n+l+s_d+j}
(R_{d-1}\vep)^{l_2+s_d-j+j_2}\sif_{d-2}^{l_2+1}
 \Zq\Big(\wvec{s_k,\dots,s_{d-3},s_{d-2}-j_2}{r_k,\,\dots\, ,r_{d-3},\, \ \ R_{d-2}\, \ };\vep\Big)
  X \rrp.  \endaligned   \label{equ:induc6}
\end{align}
For \eqref{equ:induc11} observe that $\displaystyle {\tilde P}\llp
\uZq\lp\wvec{\vs^{(p)}}{\vr^{(p)}}\rp X\rrp$ can not produce
legitimate un-tailored fractions by part (b). So the only possible
un-tailored terms are the result of negative powers of $R_{d-1}$
coming from  $\displaystyle \uZq\lp\wvec{s_{d-1}-j}{R_{d-1}}\rp$.
When $n\le 0$ the only way $R_{d-1}$ can appears in the
denominator is when $n+l+s_d+j>1-s_{d-1}$. Then these terms are
multiplied by $r_d^{l+s_d+j}$ to become tailored.

Let's consider \eqref{equ:induc5} next. Put
$w:=n+l+s_d+j-(1+j-s_{d-1}-j_2)$. We treat the two cases $w\le 0$
and $w>0$ separately. If $w\le 0$, then the induction assumption
takes care of all the terms except negative powers $R_{d-1}$. But
by breaking \eqref{equ:induc4} into two parts like in
\eqref{equ:induc0} we see that when $w\le 0$ negative powers of
$R_{d-1}$ can not appear by part (c). So we can assume now $w>0$.
As $n\le 0$ this means that $l+s_d+j\ge n+l+s_d+j>1+j-s_{d-1}-j_2$
which implies that
 $$r_d^{l+s_d+j}\vep^{n+l+s_d+j}\Big(\frac{-1}{R_{d-1}\vep}\Big)^{1+j-s_{d-1}-j_2}=
  \left(\frac{r_d}{R_{d-1}}\right)^{l+s_d+j} \vep^n (R_{d-1}\vep)^{w-n}$$
is tailored. Hence we see that \eqref{equ:induc5} is now tailored
by setting $a=w-n>0$ and $t=d-2$ in the following Claim.

\bigskip

\noindent \textbf{Claim.} Let $a$ and $t$ be two positive integers
such that $1\le t\le d-1$. Then the series
$$\sum_{\vi(t)}   \tilde{P}\llp \vep^n (R_{t+1}\vep)^a
   \Zq\Big(\wvec{s_k,\dots,s_{t-1},s_t}{r_k,\dots,r_{t-1},R_t};\vep\Big)
 X \rrp$$
is tailored.

\bigskip We use induction on $t$. Notice that this is the inner
induction loop. The outer induction loop is on $d$. When $t=1$ the
claim follows from Theorem~\ref{thm:Riemanncase} easily. This
already proves the proposition if $d=2$. We now assume $d>2$,
$t>1$, and the claim is true if the length of the vector in the
claim is less than $t$ for some $1<t<d$.

As before we may break $\Pi_t$ into two parts and form the
disjoint union
\begin{align*}
  \Pi_t=&\{(i_1,\cdots,i_p,t)\ \big|\ (i_1,\cdots,i_p)\in \Pi_{t-1}\}  \\
& \cup
    \{(i_1,\cdots,i_{p-1},i_p+1)\ \big|\ (i_1,\cdots,i_p)\in
    \Pi_{t-1}\}.
\end{align*}
Adopting the same argument we used to obtain \eqref{equ:induc0} to
\eqref{equ:induc3} we have
\begin{align}
\ &\sum_{\vi(t)}   \tilde{P}\llp \vep^n (R_{t+1}\vep)^a
   \Zq\Big(\wvec{s_k,\dots,s_{t-1},s_t}{r_k,\dots,r_{t-1},R_t};\vep\Big)
 X \rrp\\
 &=\sum_{\vi(t-1)}
  \tilde{P}\llp \vep^n  (R_{t+1}\vep)^a \uZq\lp\wvec{s_t}{R_t}\rp\check{P}\llp
\uZq\lp\wvec{\vs^{(p)}}{\vr^{(p)}}\rp X\rrp+ \vep^n
(R_{t+1}\vep)^a \uZq\lp\wvec{\vs^{(p)},s_t}{\vr^{(p)},R_t}\rp
X\rrp
\label{equ:induc01}\\
&=\sum_{\vi(t-1)}
 \tilde{P}\llp \vep^n (R_{t+1}\vep)^a\uZq\lp\wvec{s_t}{R_t}\rp{\tilde P}\llp
\uZq\lp\wvec{\vs^{(p)}}{\vr^{(p)}}\rp X\rrp \rrp\label{equ:induc7}\\
 &- \frac{q-1}{\ln q}\sum_{j=0}^{-s_t} \frac{(-s_t)!}{j!}\sum_{\vi(t-1)}
 \tilde{P}\llp \vep^n (R_{t+1}\vep)^a \Big(\frac{-1 }{R_t\vep}\Big)^{1-s_t-j}
\Zq\Big(\wvec{s_k,\dots,s_{t-2},s_{t-1}-j}{r_k,\dots,r_{t-2},\, \ \ R_{t-1}\, \ };\vep\Big)X \rrp \label{equ:induc8}\\
\ &-\sum_{j=0}^{-s_t}  {-s_t\choose j}\sum_{l=-s_t-j}^\infty
  \sum_{\vi(t-1)} \tilde{P}\llp \vep^n (R_{t+1}\vep)^a  \frac{\zq(-l)(R_t\vep)^{l+s_t+j}}{(l+s_t+j)!}
  \sif_{t-1}^{l+1}\Zq\Big(\wvec{s_k,\dots,s_{t-2},s_{t-1}-j}
  {r_k,\dots,r_{t-2},\, \ \ R_{t-1}\, \ };\vep\Big)X \rrp.
  \label{equ:induc9}
\end{align}

By part (c) $R_t$ appears in some denominator of the first term of
\eqref{equ:induc01} only if $n+a> 1-s_t$ in which case the
exponent of $R_t$ is $1-s_t$. Then multiplied by $R_{t+1}^a$ this
becomes tailored as $a\ge a+n$. Un-tailored terms with $R_t$
appearing in the denominator can not be produced by the second
term by Remark~\ref{rem:rd}.

The same argument for \eqref{equ:induc1} above using part (b) and
part (c) shows that \eqref{equ:induc7} does not contribute to any
un-tailored term in the end except possible powers of $1/R_t$.
However, all such terms are multiplied by high enough powers of
$R_{t+1}$ resulting in tailored terms again.

Similar argument for  \eqref{equ:induc2} implies that in
\eqref{equ:induc8} we only need to consider the case when $n+a>
1-s_t$. In this case we see that
 $$\vep^n (R_{t+1}\vep)^a \Big(\frac{-1 }{R_t\vep}\Big)^{1-s_t-j}
 =\vep^n \Big(\frac{-R_{t+1}}{R_t}\Big)^a\cdot
 (R_t\vep)^{a-(1-s_t-j)}$$
is tailored and therefore \eqref{equ:induc8} is tailored Laurent
series by inner loop induction assumption. The same argument by
manipulating the powers of $\vep$ works for \eqref{equ:induc9} in
the same fashion and this finishes the proof of the claim.
\medskip

Finally, \eqref{equ:induc6} follows from the claim immediately by
setting $a=l+2s_d+l_2+j_2>0$ and $t=d-1$ since
 $$\vep^{n+l+s_d+j} r_d^{l+s_d+j}
(R_{d-1}\vep)^{l_2+s_d-j+j_2}
 =\vep^n  \left(\frac{r_d}{R_{d-1}}\right)^{l+s_d+j} (R_{d-1}\vep)^{l+2s_d+l_2+j_2}. $$
This finishes the proof of the proposition.
\end{proof}

 \bigskip

\noindent \emph{Proof of Theorem \ref{thm:gzeta}}. Clearly Theorem
\ref{thm:gzeta} is true when $d=1$ by
Theorem~\ref{thm:Riemanncase}. When $d\ge 2$ by setting $m=n=0$ in
Proposition~\ref{prop:mn} in (d) we see that the constant term
$c_0$ of $\displaystyle{\tilde
Z}_{q+}^{0,0}\Big(\wvec{\vs}{\vr}\Big)$ must be the product of
linear tailored fractions of the form
$$\frac{r_i+r_{i+1}+\cdots+r_d}{r_j+r_{j+1}+\cdots+r_d}$$
for some $i$ and $j$ such that $d\ge i>j\ge 1$. Note that in the
limit \eqref{equ:shifted} each $r_k$ can only be either $\gd+n_k$,
or $2\gd+n_k$ for some nonnegative integer $n_k$. If $n_k>0$ for
some $j\le k\le d$ then the limit
 $$\lim_{\gd\to 0^+}
 \frac{r_i+r_{i+1}+\cdots+r_d}{r_j+r_{j+1}+\cdots+r_d}$$
is clearly finite. If $n_k=0$ for all $k=j,j+1,\dots,d$ then the
limit is still finite because $\gd$ will be cancelled out.

This completes the proof of Theorem \ref{thm:gzeta}.
 \hfill$\square$

 \bigskip

One of the most important properties of MZVs and M$q$ZVs is that
they satisfy the ($q$-)stuffle relations. In order to show that
our renormalization for M$q$ZVs at all non-positive arguments is
correct we need to show some type of $q$-stuffle relation hold for
them. However, we must use the shifted renormalization following
the shifting principle of M$q$ZVs.

\begin {thm}\label{thm:gshuf}
The shifted renormalizations of $\zq(\vs)$ satisfy the $q$-stuffle
relations for all $\vs\in \Z_{\leq 0}^k$. Explicitly, let $\vs_1
\in \Z_{\leq 0}^k, \vs_2 \in \Z_{\leq 0}^\ell$ and assume there
are $W_t$ vectors of length $k+\ell-t$ produced by the $q$-stuffle
and let $\vf_j$ be the binary vector representing the positions of
the stuffing in $j$-th such vector:
$$ \vs_1\msh \vs_2=\sum_{t=0}^{\min(k,\ell)}
 (1-q)^t\sum_{j=1}^{W_t} (\vec{v}_j-\vf_j).$$
Then
$$\gzq(\vs_1)\gzq(\vs_2)=\sum_{t=0}^{\min{k,\ell}}
 (1-q)^t\sum_{j=1}^{W_t} \gzq^{\vf_j}(\vec{v}_j-\vf_j).$$
When $\qup$ we obtain \cite[Theorem~4.11]{GuZ}.
\end{thm}
\begin{proof} By the definition of renormalization in
Definition \ref{defn:noshiftMqZV} and Corollary~\ref{cor:stuffle}
we have
 \begin{equation*}
 \gzq(\vs_1)\gzq(\vs_2)=\lim_{\vr_1 \to -\vs_1^+, \vr_2\to -\vs_2^+}
\zeta\lp\wvec{\vs_1}{\vr_1}\rp \zeta\lp\wvec{\vs_2}{\vr_2}\rp
=\lim_{\vr_1 \to -\vs_1^+, \vr_2\to -\vs_2^+}
\zeta\lp\wvec{\vs_1}{\vr_1}\msh\wvec{\vs_2}{\vr_2}\rp.
\end{equation*}
 From the definition of $q$-stuffle \eqref{equ:stuf} and the shifted
renormalization \eqref{equ:shifted} the theorem follows
immediately.
\end{proof}

To conclude our paper we observe that a multiple $q$-zeta function
has more singularities than its classical counterpart. We don't
know how to renormalize M$q$ZVs at these points at present.

\noindent {\em Email:} zhaoj@eckerd.edu

\end{document}